\documentclass[12pt]{article}
\usepackage{amsfonts}
\usepackage{amssymb}
\usepackage{graphicx}

\begin{document}

\author{Francis OGER}
\title{\textbf{Equivalence \'{e}l\'{e}mentaire entre pavages}}
\date{}

\begin{center}
\textbf{Paperfolding sequences, paperfolding curves}

\textbf{and local isomorphism}

\bigskip

Francis OGER
\end{center}

\bigskip

\bigskip

\noindent ABSTRACT. For each integer $n$, an $n$-folding curve is
obtained by folding $n$ times a strip of paper in two, possibly up or down,
and unfolding it with right angles. Generalizing the usual notion of
infinite folding curve, we define complete folding curves as the curves
without endpoint which are unions of increasing sequences of $n$-folding
curves for $n$ integer.

We prove that there exists a standard way to extend any complete folding
curve into a covering of $%
\mathbb{R}
^{2}$ by disjoint such curves, which satisfies the local isomorphism
property introduced to investigate aperiodic tiling systems. This covering
contains at most six curves.\bigskip 


\bigskip

\bigskip

The infinite folding sequences (resp. curves) usually considered are
sequences $(a_{k})_{k\in 
\mathbb{N}
^{\ast }}\subset \{+1,-1\}$ (resp. infinite curves with one endpoint)
obtained as direct limits of $n$-folding sequences (resp. curves) for $n\in 
\mathbb{N}
$. It is well known (see [3] and [4]) that paperfolding curves are
self-avoiding and that, in some cases, including the Heighway Dragon curve,
a small number of copies of the same infinite folding curve can be used to
cover $%
\mathbb{R}
^{2}$ without overlapping. Anyway, in some other cases, the last property is
not true.

In the present paper, we define complete folding sequences (resp. curves) as
the sequences $(a_{k})_{k\in 
\mathbb{Z}
}\subset \{+1,-1\}$ (resp. the infinite curves without endpoint) which are
direct limits of $n$-folding sequences (resp. curves) for $n\in 
\mathbb{N}
$. Any infinite folding sequence (resp. curve) in the classical sense can be
extended into a complete folding sequence (resp. curve). On the other hand,
most of the complete folding sequences (resp. curves) cannot be obtained in
that way.

We prove that any complete folding curve, and therefore any infinite folding
curve, can be extended in an essentially unique way into a covering of $%
\mathbb{R}
^{2}$ by disjoint complete folding curves which satisfies the local
isomorphism property. We show that a covering obtained from an infinite
folding curve can contain complete folding curves which are not extensions
of infinite folding curves.

One important argument in the proofs is the derivation of paperfolding
curves, which is investigated in Section 2. Another one is the local
isomorphism property for complete folding sequences (cf. Section 1) and for
coverings of $%
\mathbb{R}
^{2}$ by sets of disjoint complete folding curves (cf. Section 3). The local
isomorphism property was originally used to investigate aperiodic tiling
systems. Actually, we have an interpretation of complete folding sequences
as tilings of $%
\mathbb{R}
$, and an interpretation of coverings of $%
\mathbb{R}
^{2}$ by disjoint complete folding curves as tilings of $%
\mathbb{R}
^{2}$.\bigskip

\textbf{1. Paperfolding sequences.}\bigskip

The notions usually considered (see for instance [5]), and which we define
first, are those of $n$-folding sequence (sequence obtained by folding $n$
times a strip of paper in two), and $\infty $-folding sequence (sequence
indexed by $%
\mathbb{N}
^{\ast }$ which is obtained as a direct limit of $n$-folding sequences for $%
n\in 
\mathbb{N}
$).

Then we introduce complete folding sequences, which are sequences indexed by 
$%
\mathbb{Z}
$, also obtained as direct limits of $n$-folding sequences for $n\in 
\mathbb{N}
$. We describe the finite subwords of each such sequence. Using this
description, we show that complete folding sequences satisfy properties
similar to those of aperiodic tiling systems: they form a class defined by a
set of local rules, neither of them is periodic, but all of them satisfy the
local isomorphism property introduced for tilings. It follows that, for each
such sequence, there exist $2^{\omega }$ isomorphism classes of sequences
which are locally isomorphic to it.\bigskip

\noindent \textbf{Definitions.} For each $n\in 
\mathbb{N}
$\ and each sequence $S=(a_{1},...,a_{n})\subset \{+1,-1\}$, we write $%
\left\vert S\right\vert =n$\ and $\overline{S}=(-a_{n},...,-a_{1})$. We say
that a sequence $(a_{1},...,a_{n})$ is a \emph{subword} of a sequence $%
(b_{1},...,b_{p})$ or $(b_{k})_{k\in 
\mathbb{N}
^{\ast }}$ or $(b_{k})_{k\in 
\mathbb{Z}
}$ if there exists $h$ such that $a_{k}=b_{k+h}$ for $1\leq k\leq n$.\bigskip

\noindent \textbf{Definition.} For each $n\in 
\mathbb{N}
$, an $n$\emph{-folding sequence} is a sequence

\noindent $(a_{1},...,a_{2^{n}-1})\subset \{+1,-1\}$\ obtained by folding $n$%
\ times a strip of paper in two, with each folding being done independently
up or down, unfolding it, and writing $a_{k}=+1$\ (resp. $a_{k}=-1$) for
each $k\in \{1,...,2^{n}-1\}$\ such that the $k$-th fold from the left has
the shape of a $\vee $\ (resp. $\wedge $) (we obtain\ the empty sequence for 
$n=0$).\bigskip

\noindent \textbf{Properties.} The following properties are true for each $%
n\in 
\mathbb{N}
$:

\noindent 1) If $S$ is an $n$-folding sequence, then $\overline{S}$ is also
an $n$-folding sequence.

\noindent 2) The $(n+1)$-folding sequences are the sequences $(\overline{S}%
,+1,S)$ and the sequences $(\overline{S},-1,S)$, where $S$ is an $n$-folding
sequence.

\noindent 3) There exist $2^{n}$ $n$-folding sequences (proof by induction
on $n$ using 2)).

\noindent 4) Any sequence $(a_{1},...,a_{2^{n+1}-1})$ is a $(n+1)$-folding
sequence if and only if $(a_{2k})_{1\leq k\leq 2^{n}-1}$\ is an $n$-folding
sequence and $a_{1+2k}=(-1)^{k}a_{1}$ for $0\leq k\leq 2^{n}-1$.

\noindent 5) If $n\geq 2$ and if $(a_{1},...,a_{2^{n}-1})$ is an $n$-folding
sequence, then $a_{2^{r}(1+2k)}=(-1)^{k}a_{2^{r}}$ for $0\leq r\leq n-2$ and 
$0\leq k\leq 2^{n-r-1}-1$ (proof by induction on $n$ using 4)).\bigskip

\noindent \textbf{Definition.} An $\infty $\emph{-folding sequence} is a
sequence $(a_{n})_{n\in 
\mathbb{N}
^{\ast }}$ such that

\noindent $(a_{1},...,a_{2^{n}-1})$ is an $n$-folding sequence for each $%
n\in 
\mathbb{N}
^{\ast }$.\bigskip

\noindent \textbf{Definition.} A \emph{finite folding sequence} is a subword
of an $n$-folding sequence for an integer $n$.\bigskip

\noindent \textbf{Examples.} The sequence $(+1,+1,+1)$\ is a finite folding
sequence since it is a subword of the $3$-folding sequence $%
(-1,+1,+1,+1,-1,-1,+1)$. On the other hand, $(+1,+1,+1)$ is not a $2$%
-folding sequence and $(+1,+1,+1,+1)$\ is not a finite folding
sequence.\bigskip

\noindent \textbf{Definition.} A\emph{\ complete folding sequence} is a
sequence $(a_{k})_{k\in 
\mathbb{Z}
}\subset \{+1,-1\}$ such that its finite subwords are folding
sequences.\bigskip

\noindent \textbf{Examples.} For each $\infty $-folding sequence $%
S=(a_{n})_{n\in 
\mathbb{N}
^{\ast }}$, write $\overline{S}=(-a_{-n})_{n\in -%
\mathbb{N}
^{\ast }}$. Then $(\overline{S},+1,S)$ and $(\overline{S},-1,S)$ are
complete folding sequences since

\noindent $(-a_{2^{n}-1},...,-a_{1},+1,a_{1},...,a_{2^{n}-1})$ and $%
(-a_{2^{n}-1},...,-a_{1},-1,a_{1},...,a_{2^{n}-1})$

\noindent are $(n+1)$-folding sequences for each $n\in 
\mathbb{N}
$. In Section 3, we give examples of complete folding sequences which are
not obtained in that way.\bigskip

It follows from the property 5) above that, for each complete folding
sequence $(a_{h})_{h\in 
\mathbb{Z}
}$ and each $n\in 
\mathbb{N}
$, there exists $k\in 
\mathbb{Z}
$\ such that $a_{k+l.2^{n+1}}=(-1)^{l}a_{k}$\ for each $l\in 
\mathbb{Z}
$. Moreover we have:\bigskip

\noindent \textbf{Proposition 1.1.} Consider a sequence $S=(a_{h})_{h\in 
\mathbb{Z}
}\subset \{+1,-1\}$. For each $n\in 
\mathbb{N}
$, suppose that there exists $h_{n}\in 
\mathbb{Z}
$\ such that $a_{h_{n}+k.2^{n+1}}=(-1)^{k}a_{h_{n}}$\ for each $k\in 
\mathbb{Z}
$, and consider $E_{n}=h_{n}+2^{n}%
\mathbb{Z}
$ and $F_{n}=h_{n}+2^{n+1}%
\mathbb{Z}
$.\ Then, for each $n\in 
\mathbb{N}
$:

\noindent 1) $%
\mathbb{Z}
-E_{n}=\{h\in 
\mathbb{Z}
$~$\mathbf{\mid }$ $a_{h+k.2^{n+1}}=a_{h}$ for each $k\in 
\mathbb{Z}
\mathbf{\}}$ and $%
\mathbb{Z}
-E_{n}$ is the disjoint union of $F_{0},...,F_{n-1}$;

\noindent 2) for each $h\in E_{n}$, $(a_{h-2^{n}+1},...,a_{h+2^{n}-1})$\ is
a $(n+1)$-folding sequence;

\noindent 3) for each $h\in 
\mathbb{Z}
$, if $(a_{h-2^{n+1}+1},...,a_{h+2^{n+1}-1})$\ is a $(n+2)$-folding
sequence, then $h\in E_{n}$.\bigskip

\noindent \textbf{Proof.} It follows from the definition of the integers $%
h_{n}$\ that the sets $F_{n}$\ are disjoint. For each $n\in 
\mathbb{N}
$, we have $E_{n}=%
\mathbb{Z}
-(F_{0}\cup ...\cup F_{n-1})$ since $%
\mathbb{Z}
-(F_{0}\cup ...\cup F_{n-1})$ is of the form $h+2^{n}%
\mathbb{Z}
$ and $h_{n}$\ does not belong to $F_{0}\cup ...\cup F_{n-1}$.

For each $n\in 
\mathbb{N}
$, there exists no $h\in E_{n}$\ such that $a_{h+k.2^{n+1}}=a_{h}$ for each $%
k\in 
\mathbb{Z}
$, since we have $E_{n}=(h_{n}+2^{n+1}%
\mathbb{Z}
)\cup (h_{n+1}+2^{n+1}%
\mathbb{Z}
)$, $a_{h_{n}+2^{n+1}}=-a_{h_{n}}$ and $a_{h_{n+1}+2^{n+2}}=-a_{h_{n+1}}$.
On the other hand, we have $a_{h_{m}+k.2^{m+2}}=a_{h_{m}}$\ for $0\leq m\leq
n-1$ and $k\in 
\mathbb{Z}
$, and therefore $a_{h+k.2^{n+1}}=a_{h}$ for $h\in F_{0}\cup ...\cup F_{n-1}$
and $k\in 
\mathbb{Z}
$, which completes the proof de 1).

We show 2) by induction on $n$. The case $n=0$\ is clear. If 2) is true for $%
n$, then, for each $h\in E_{n+1}$, the induction hypothesis applied to $%
(a_{h+2k})_{k\in 
\mathbb{Z}
}$ implies that $(a_{h-2^{n+1}+2k})_{1\leq k\leq 2^{n+1}-1}$\ is a $(n+1)$%
-folding sequence; it follows that $(a_{h-2^{n+1}+1},...,a_{h+2^{n+1}-1})$\
is a $(n+2)$-folding sequence, since $%
a_{h-2^{n+1}+1+2k}=(-1)^{k}a_{h-2^{n+1}+1}$\ for $0\leq k\leq 2^{n+1}-1$.

Concerning 3), we observe that, for each $h\in 
\mathbb{Z}
$, if $(a_{h-2^{n+1}+1},...,a_{h+2^{n+1}-1})$\ is a $(n+2)$-folding
sequence, then $a_{h+2^{n}}=-a_{h-2^{n}}$. According to 1), it follows $%
h-2^{n}\in E_{n}$, and therefore\ $h\in E_{n}$.~~$\blacksquare $\bigskip

\noindent \textbf{Corollary 1.2.} Any sequence $(a_{h})_{h\in 
\mathbb{Z}
}\subset \{+1,-1\}$ is a complete folding sequence if and only if, for each $%
n\in 
\mathbb{N}
$, there exists $h_{n}\in 
\mathbb{Z}
$\ such that $a_{h_{n}+k.2^{n+1}}=(-1)^{k}a_{h_{n}}$\ for each $k\in 
\mathbb{Z}
$.\bigskip

For each complete folding sequence $S$ and each $n\in 
\mathbb{N}
$, the sets $E_{n}$ and $F_{n}$\ of Proposition 1.1 do not depend on the
choice of $h_{n}$. We denote them by $E_{n}(S)$ and $F_{n}(S)$.\ We write $%
E_{n}$ and $F_{n}$\ instead of $E_{n}(S)$ and $F_{n}(S)$ if it creates no
ambiguity.\bigskip

\noindent \textbf{Corollary 1.3.} Any complete folding sequence is
nonperiodic.\bigskip

\noindent \textbf{Proof. }Let $S=(a_{h})_{h\in 
\mathbb{Z}
}$ be such a sequence, and let $r$ be an integer such that $a_{h+r}=a_{h}$\
for each $h\in 
\mathbb{Z}
$.\ For each $n\in 
\mathbb{N}
$, it follows from 1) of Proposition 1.1 that $r+(%
\mathbb{Z}
-E_{n})=%
\mathbb{Z}
-E_{n}$, whence $r+E_{n}=E_{n}$ and $r\in 2^{n}%
\mathbb{Z}
$. Consequently, we have $r=0$.~~$\blacksquare $\bigskip

Now, for each complete folding sequence $S=(a_{h})_{h\in 
\mathbb{Z}
}$, we describe the finite subwords of $S$ and we count those which have a
given length.\bigskip

\noindent \textbf{Lemma 1.4.} For each $n\in 
\mathbb{N}
$ and for any $r,s\in 
\mathbb{Z}
$, we have $r-s\in 2^{n+1}%
\mathbb{Z}
$ if $(a_{r+1},...,a_{r+t})=(a_{s+1},...,a_{s+t})$ for $t=\sup (2^{n},7)$%
.\bigskip

\noindent \textbf{Proof.} If $r-s\notin 2%
\mathbb{Z}
$, then we have for instance $r\in E_{1}$ and $s\in F_{0}$. It follows $%
a_{r+1}=-a_{r+3}=a_{r+5}=-a_{r+7}$\ since $r+1\in F_{0}$. Moreover, we have $%
a_{s+5}=-a_{s+1}$ if $s+1\in F_{1}$, and $a_{s+7}=-a_{s+3}$ if $s+3\in F_{1}$%
. One of these two possibilities is necessarily realized\ since $s+1\in
E_{1} $, which contradicts $(a_{r+1},...,a_{r+7})=(a_{s+1},...,a_{s+7})$.

If $r-s\in 2^{k}%
\mathbb{Z}
-2^{k+1}%
\mathbb{Z}
$ with $1\leq k\leq n$, then we consider $h\in \{1,...,2^{k}\}$ such that $%
r+h\in F_{k-1}$. We have $a_{r+h+m.2^{k}}=(-1)^{m}a_{r+h}$ for each $m\in 
\mathbb{Z}
$, and in particular $a_{s+h}=-a_{r+h}$, which contradicts $%
(a_{r+1},...,a_{r+t})=(a_{s+1},...,a_{s+t})$\ since $1\leq h\leq 2^{k}\leq t$%
.~~$\blacksquare $\bigskip

\noindent \textbf{Proposition 1.5.} Consider $n\in 
\mathbb{N}
$ and write\ $T=(a_{h+1},...,a_{h+2^{n}-1})$ with $h\in E_{n}$. Then any
sequence of length $\leq 2^{n+1}-1$\ is a subword de $S$ if and only if
there exist $\zeta ,\eta \in \{-1,+1\}$ such that it can be written in one
of the forms:

\noindent (1) $(T_{1},\zeta ,T_{2})$ with $T_{1}$ final segment of $T$ and $%
T_{2}$ initial segment of $\overline{T}$;

\noindent (2) $(T_{1},\zeta ,T_{2})$ with $T_{1}$ final segment of $%
\overline{T}$ and $T_{2}$ initial segment of $T$;

\noindent (3) $(T_{1},\zeta ,\overline{T},\eta ,T_{2})$ with $T_{1}$ final
segment and $T_{2}$ initial segment of $T$;

\noindent (4) $(T_{1},\zeta ,T,\eta ,T_{2})$ with $T_{1}$ final segment and $%
T_{2}$ initial segment of $\overline{T}$.

\noindent If $\sup (2^{n},7)\leq t\leq 2^{n+1}-1$, then any subword of
length $t$ of $S$\ can be written in exactly one way in one of the forms
(1), (2), (3), (4).\bigskip

\noindent \textbf{Proof.} We can suppose $h\in F_{n}$ since $T$ and $%
\overline{T}$\ play symmetric roles in the Proposition. Then we have $%
(a_{k+1},...,a_{k+2^{n}-1})=T$ for $k\in F_{n}$ and $%
(a_{k+1},...,a_{k+2^{n}-1})=\overline{T}$ for $k\in E_{n+1}$. It follows
that each subword of $S$ of length $\leq 2^{n+1}-1$ can be written in one of
the forms (1), (2), (3), (4).

Now, we are going to prove that each sequence of one of these forms can be
expressed as a subword of $S$ in such a way that the part $T_{1}$ is
associated to the final segment of a sequence $(a_{k+1},...,a_{k+2^{n}-1})$
with $k\in E_{n}$.

First we show this property for the sequences of the form (1) or (3). It
suffices to prove that, for any $\zeta ,\eta \in \{-1,+1\}$, there exists $%
k\in F_{n}$ such that $a_{k-2^{n}}=\zeta $ and $a_{k}=\eta $, since these
two equalities imply $(a_{k-2^{n+1}+1},...,a_{k+2^{n}-1})=(T,\zeta ,%
\overline{T},\eta ,T)$. We consider $l\in F_{n}$\ such that $a_{l}=\eta $.
We have $a_{l+r.2^{n+2}}=\eta $\ for each $r\in 
\mathbb{Z}
$. Moreover, $\{l+r.2^{n+2}-2^{n}\mid r\in 
\mathbb{Z}
\}$ is equal to $F_{n+1}$ or $E_{n+2}$. In both cases, there exists $r\in 
\mathbb{Z}
$ such that $a_{l+r.2^{n+2}-2^{n}}=\zeta $, and it suffices to take $%
k=l+r.2^{n+2}$ for such an $r$.

Now we show the same property for the sequences of the form (2) or (4). It
suffices to prove that, for any $\zeta ,\eta \in \{-1,+1\}$, there exists $%
k\in F_{n}$ such that $a_{k}=\zeta $ and $a_{k+2^{n}}=\eta $, since these
two equalities imply $(a_{k-2^{n}+1},...,a_{k+2^{n+1}-1})=(\overline{T}%
,\zeta ,T,\eta ,\overline{T})$. We consider $l\in F_{n}$\ such that $%
a_{l}=\zeta $. We have $a_{l+r.2^{n+2}}=\zeta $\ for each $r\in 
\mathbb{Z}
$. Moreover, $\{l+r.2^{n+2}+2^{n}\mid r\in 
\mathbb{Z}
\}$ is equal to $F_{n+1}$ or $E_{n+2}$. In both cases, there exists $r\in 
\mathbb{Z}
$ such that $a_{l+r.2^{n+2}+2^{n}}=\eta $, and it suffices to take $%
k=l+r.2^{n+2}$ for such an $r$.

Now, suppose that two expressions of the forms (1), (2), (3), (4) give the
same sequence of length $t$ with $\sup (2^{n},7)\leq t\leq 2^{n+1}-1$.
Consider two sequences $(a_{r+1},...,a_{r+t})$ and $(a_{s+1},...,a_{s+t})$
which realize these expressions in such a way that, in each of them, the
part $T_{1}$ of the expression is associated to a final segment of a
sequence $(a_{k+1},...,a_{k+2^{n}-1})$ with $k\in E_{n}$, while the part $%
T_{2}$ is associated to an initial segment of a sequence $%
(a_{l+1},...,a_{l+2^{n}-1})$ with $l=k+2^{n}$ or $l=k+2^{n+1}$. Then, by\
Lemma 1.4, the equality $(a_{r+1},...,a_{r+t})=(a_{s+1},...,a_{s+t})$\
implies $r-s\in 2^{n+1}%
\mathbb{Z}
$.\ It follows that the two expressions are equal.~~$\blacksquare $\bigskip

\noindent \textbf{Corollary 1.6.} Any finite folding sequence $U$ is a
subword of $S$ if and only if $\overline{U}$ is a subword of $S$.\bigskip

\noindent \textbf{Proof.} For each sequence $T=(a_{h+1},...,a_{h+2^{n}-1})$
with $n\in 
\mathbb{N}
$ and $h\in E_{n}$, the sequence $U$ is of the form (1) (resp. (2), (3),
(4)) relative to $T$ if and only if $\overline{U}$ is of the form (2) (resp.
(1), (4), (3)) relative to $T$.~~$\blacksquare $\bigskip

It follows from the Corollary below that, for each integer $n\geq 3$, each
complete folding sequence has exactly $8$ subwords which are $n$-folding
sequences:\bigskip

\noindent \textbf{Corollary 1.7.} Consider $n\in 
\mathbb{N}
$ and write\ $T=(a_{h+1},...,a_{h+2^{n}-1})$ with $h\in E_{n}$. Then any $%
(n+2)$-folding sequence is a subword of $S$ if and only if it can be written
in the form $(\overline{T},-\zeta ,T,\eta ,\overline{T},\zeta ,T)$ or $%
(T,-\zeta ,\overline{T},\eta ,T,\zeta ,\overline{T})$ with $\zeta ,\eta \in
\{+1,-1\}$.\bigskip

\noindent \textbf{Proof. }For each $k\in 
\mathbb{Z}
$, if $(a_{k-2^{n+1}+1},...,a_{k+2^{n+1}-1})$\ is a $(n+2)$-folding
sequence, then $k\in E_{n}$ by 3) of Proposition 1.1. Consequently, we have $%
(a_{k+1},...,a_{k+2^{n}-1})=T$ or $(a_{k+1},...,a_{k+2^{n}-1})=\overline{T}$%
, and $(a_{k-2^{n+1}+1},...,a_{k+2^{n+1}-1})$ is of the required form.

In order to prove that each sequence of that form is a subword of $S$, we
consider $k\in E_{n+1}$ and we write $U=(a_{k+1},...,a_{k+2^{n+1}-1})$. We
have $U=(\overline{T},\varepsilon ,T)$ or $U=(T,\varepsilon ,\overline{T})$
with $\varepsilon =\mp 1$. Here we only consider the first case; the second
one can be treated in the same way since $T$ and $\overline{T}$ play
symmetric roles in the Corollary.

We apply Proposition 1.5 for $n+1$ instead of $n$, and we consider the forms
(1), (2), (3), (4) relative to $U$. For any $\zeta ,\eta \in \{+1,-1\}$, the
sequence $(\overline{T},-\zeta ,T,\eta ,\overline{T},\zeta ,T)$ is a subword
of $S$ because it is equal to $(U,\eta ,\overline{U})$ or to $(\overline{U}%
,\eta ,U)$, and therefore of the form (1) or (2) relative to $U$. The
sequence $(T,-\zeta ,\overline{T},\eta ,T,\zeta ,\overline{T})$ is also a
subword of $S$ because it is of the form $(T,\alpha ,\overline{U},\beta ,%
\overline{T})$ or $(T,\alpha ,U,\beta ,\overline{T})$ with $\alpha ,\beta
\in \{+1,-1\}$, and therefore of the form (3) or (4) relative to $U$.~~$%
\blacksquare $\bigskip

The following result generalizes [1, Th., p. 27] to complete folding
sequences:\bigskip

\noindent \textbf{Theorem 1.8.} The sequence $S$ has $4t$ subwords of length 
$t$ for each integer $t\geq 7$ and $2$, $4$, $8$, $12$, $18$, $23$ subwords
of length $t=1$, $2$, $3$, $4$, $5$, $6$.\bigskip

\noindent \textbf{Proof.} The proof of the Theorem for $t=1$, $2$, $3$, $4$, 
$5$, $6$ is based on Proposition 1.5. We leave it to the reader.

For $t\geq 7$, we consider the integer $n\geq 2$ such that $2^{n}\leq t\leq
2^{n+1}-1$, and we write $T=(a_{h+1},...,a_{h+2^{n}-1})$ with $h\in E_{n}$.\
By Proposition 1.5, it suffices to count the subwords of length $t$ of $S$
which are in each of the forms (1), (2), (3), (4) relative to $T$.

Each of the forms (1), (2) gives $\left\vert T_{1}\right\vert +\left\vert
T_{2}\right\vert =t-1$, and therefore $\left\vert T_{1}\right\vert \geq
(t-1)-(2^{n}-1)=t-2^{n}$. As $\left\vert T_{1}\right\vert \leq 2^{n}-1$, we
have $(2^{n}-1)-(t-2^{n})+1=2^{n+1}-t$\ possible values for $\left\vert
T_{1}\right\vert $. Consequently, there exist $4(2^{n+1}-t)$\ sequences
associated to these two forms, since there are $2$ possible values for $%
\zeta $.

Each of the forms (3), (4) gives $\left\vert T_{1}\right\vert +\left\vert
T_{2}\right\vert =t-2^{n}-1$, and therefore $\left\vert T_{1}\right\vert
\leq t-2^{n}-1$.\ We have $t-2^{n}$\ possible values for $\left\vert
T_{1}\right\vert $. Consequently, there exist $8(t-2^{n})$\ sequences
associated to these two forms, since there are $4$ possible values for $%
(\zeta ,\eta )$. Now, the total number of subwords of length $t$ in $S$ is $%
4(2^{n+1}-t)+8(t-2^{n})=4t$.~~$\blacksquare $\bigskip

For each sequence $(a_{h})_{h\in 
\mathbb{Z}
}\subset \{+1,-1\}$, we define a tiling of $%
\mathbb{R}
$ as follows: the tiles are the intervals $[k,k+1]$ for $k\in 
\mathbb{Z}
$, where the \textquotedblleft colour\textquotedblright\ of the endpoint $k$%
\ (resp. $k+1$) is\ the sign of $a_{k}$\ (resp. $a_{k+1}$). Each tile is of
one of the forms $[+,+]$, $[+,-]$, $[-,+]$, $[-,-]$, and each pair of
consecutive tiles is of one of the forms $([+,+],[+,+])$, $([+,+],[+,-])$, $%
([+,-],[-,+])$, $([+,-],[-,-])$, $([-,+],[+,+])$, $([-,+],[+,-])$, $%
([-,-],[-,+])$, $([-,-],[-,-])$.

Concerning the theory of tilings, the reader is referred to [7], which
presents classical results and gives generalizations based on mathematical
logic. Two tilings of $%
\mathbb{R}
^{n}$ are said to be \emph{isomorphic} if they are equivalent up to
translation, and \emph{locally isomorphic} if they contain the same bounded
sets of tiles modulo translations.

We say that two sequences $(a_{h})_{h\in 
\mathbb{Z}
},(b_{h})_{h\in 
\mathbb{Z}
}\subset \{+1,-1\}$ are \emph{isomorphic} (resp. \emph{locally isomorphic})
if they are equivalent up to translation (resp. they have the same finite
subwords). This property is true if and only if the associated tilings are
isomorphic (resp. locally isomorphic). It follows from the definitions that
any sequence $(a_{h})_{h\in 
\mathbb{Z}
}\subset \{+1,-1\}$ is a complete folding sequence if it is locally
isomorphic to such a sequence.\bigskip

\noindent \textbf{Corollary 1.9.} Any complete folding sequence $%
S=(a_{h})_{h\in 
\mathbb{Z}
}$ is locally isomorphic to $\overline{S}=(-a_{-h})_{h\in 
\mathbb{Z}
}$, but not locally isomorphic to $-S=(-a_{h})_{h\in 
\mathbb{Z}
}$.\bigskip

\noindent \textbf{Proof.} The first statement is a consequence of Corollary
1.6 since each finite folding sequence $T$\ is a subword of $S$ if and only
if $\overline{T}$\ is a subword of $\overline{S}$.

In order to prove the second statement, we consider $%
T=(a_{h+1},a_{h+2},a_{h+3})$ with $h\in E_{2}$. We have $-T\neq \overline{T}$%
\ since $a_{h+3}=-a_{h+1}$. By Corollary 1.7, the $4$-folding sequences
which are subwords of $S$ are the sequences $(\overline{T},-\zeta ,T,\eta ,%
\overline{T},\zeta ,T)$ and $(T,-\zeta ,\overline{T},\eta ,T,\zeta ,%
\overline{T})$ for $\zeta ,\eta \in \{+1,-1\}$, while the $4$-folding
sequences which are subwords of $-S$ are the sequences $(-\overline{T}%
,-\zeta ,-T,\eta ,-\overline{T},\zeta ,-T)$ and $(-T,-\zeta ,-\overline{T}%
,\eta ,-T,\zeta ,-\overline{T})$ for $\zeta ,\eta \in \{+1,-1\}$.\
Consequently, $S$ and $-S$ have no $4$-folding sequence in common.~~$%
\blacksquare $\bigskip

\noindent \textbf{Remark.} For each $\infty $-folding sequence $S$, it
follows from Corollary 1.9 that $T=(\overline{S},+1,S)$ and $U=(\overline{S}%
,-1,S)$ are locally isomorphic, since $U=\overline{T}$.\bigskip

We say that a tiling $\mathcal{T}$ of $%
\mathbb{R}
^{n}$ \emph{satisfies the local isomorphism property} if, for each bounded
set of tiles $\mathcal{F}\subset \mathcal{T}$, there exists $r\in 
\mathbb{R}
_{\ast }^{+}$ such that each ball of radius $r$ in $%
\mathbb{R}
^{n}$ contains the image of $\mathcal{F}$ under a translation. Then any
tiling $\mathcal{U}$ is locally isomorphic to $\mathcal{T}$ if each bounded
set of tiles contained in $\mathcal{U}$ is the image under a translation of
a set of tiles contained in $\mathcal{T}$.

We say that a sequence $(a_{h})_{h\in 
\mathbb{Z}
}\subset \{+1,-1\}$ \emph{satisfies the local isomorphism property} if the
associated tiling satisfies the local isomorphism property.

Like Robinson tilings and Penrose tilings, complete folding sequences are 
\emph{aperiodic} in the following sense:

\noindent 1) they form a class defined by a set of rules which can be
expressed by first-order sentences (for each $n\in 
\mathbb{N}
$, we write a sentence which says that each subword of length $2^{n}$ of the
sequence considered is a subword of a $(n+1)$-folding sequence);

\noindent 2) neither of them is periodic, but all of them satisfy the local
isomorphism property.

\noindent The second statement of 2) follows from the Theorem below:\bigskip

\noindent \textbf{Theorem 1.10.} Let $S=(a_{h})_{h\in 
\mathbb{Z}
}$\ be a complete folding sequence, let $T$ be a finite subword of $S$, and
let $r$ be an integer such that $\left\vert T\right\vert \leq 2^{r}$.\ Then $%
T$ is a subword of $(a_{h+1},...,a_{h+10.2^{r}-2})$ for each $h\in 
\mathbb{Z}
$.\bigskip

\noindent \textbf{Proof. }There exists $k\in E_{r}$ such that $T$ is a
subword of $(a_{k-2^{r}+1},...,a_{k+2^{r}-1})$. We have $%
(a_{k-2^{r}+1},...,a_{k+2^{r}-1})=(\overline{U},\zeta ,U)$ with $\zeta
=a_{k} $ and $U=(a_{k+1},...,a_{k+2^{r}-1})$.

If $k\in E_{r+1}$, we consider $m\in F_{r+1}$ such that $a_{m}=\zeta $; we
have $a_{m+n.2^{r+2}}=(-1)^{n}\zeta $\ for each $n\in 
\mathbb{Z}
$. If $k\in F_{r}$, we write $m=k$; we have $a_{m+n.2^{r+1}}=(-1)^{n}\zeta $%
\ for each $n\in 
\mathbb{Z}
$. In both cases, for each $n\in 
\mathbb{Z}
$, we have $a_{m+n.2^{r+3}}=\zeta $ and $%
(a_{m+n.2^{r+3}-2^{r}+1},...,a_{m+n.2^{r+3}+2^{r}-1})=(\overline{U},\zeta
,U) $.

For each $h\in 
\mathbb{Z}
$, there exists $n\in 
\mathbb{Z}
$ such that $h-m+2^{r}\leq n.2^{r+3}\leq h-m+9.2^{r}-1$, which implies $%
h+1\leq m+n.2^{r+3}-2^{r}+1$\ and $h+10.2^{r}-2\geq m+n.2^{r+3}+2^{r}-1$.
Then $(a_{m+n.2^{r+3}-2^{r}+1},...,a_{m+n.2^{r+3}+2^{r}-1})$\ is a subword
of $(a_{h+1},...,a_{h+10.2^{r}-2})$,\ which completes the proof of the
Theorem since $T$ is a subword of $%
(a_{m+n.2^{r+3}-2^{r}+1},...,a_{m+n.2^{r+3}+2^{r}-1})$.~~$\blacksquare $%
\bigskip

The second part of the Theorem below is similar to results which were proved
for Robinson tilings and Penrose tilings:\bigskip

\noindent \textbf{Theorem 1.11.} 1) There exist $2^{\omega }$ complete
folding sequences which are pairwise not locally isomorphic.

\noindent 2) For each complete folding sequence $S$, there exist $2^{\omega
} $ isomorphism classes of sequences which are locally isomorphic to $S$%
.\bigskip

\noindent \textbf{Proof of 1).} It follows from Proposition 1.1 that each
complete folding sequence $(a_{h})_{h\in 
\mathbb{Z}
}$ is completely determined by the following operations:

\noindent - successively for each $n\in 
\mathbb{N}
$, we choose among the $2$ possible values the smallest $h\in 
\mathbb{N}
\cap F_{n}$, then we fix $a_{h}\in \{+1,-1\}$;

\noindent - for the unique $h\in \cap _{n\in 
\mathbb{N}
}E_{n}$\ if it exists, we fix $a_{h}\in \{+1,-1\}$.

\noindent Moreover, each possible sequence of choices determines a complete
folding sequence.

Now, it follows from Corollary 1.7 that, for each complete folding sequence $%
S$ and each integer $m$, there exist an integer $n>m$ and a complete folding
sequence $T$ such that $S$ and $T$ contain as subwords the same $m$-folding
sequences, but not the same $n$-folding sequences.

\noindent \textbf{Proof of 2).} The sequence $S$ is not periodic by
Corollary 1.3, and satisfies the local isomorphism property according to
Theorem 1.10. By [7, Corollary 3.7], it follows that there exist $2^{\omega
} $ isomorphism classes of sequences which are locally isomorphic to $S$.~~$%
\blacksquare $\bigskip

\noindent \textbf{Remark.} Concerning logic, we note two differences between
complete folding sequences and Robinson or Penrose tilings. First, the set
of all complete folding sequences is defined by a countable set of
first-order sentences, and not by only one sentence. Second, it is the union
of $2^{\omega }$ classes for elementary equivalence, i.e. local isomorphism,
instead of being a single class.\bigskip

\textbf{2. Paperfolding curves: self-avoiding, derivatives, exterior.}%
\bigskip

In the present section, we define $n$-folding curves, finite folding curves, 
$\infty $-folding curves and complete folding curves associated to $n$%
-folding sequences, finite folding sequences, $\infty $-folding sequences
and complete folding sequences. We show their classical properties:
self-avoiding, existence of \textquotedblleft derivatives\textquotedblright .

Then we prove that any complete folding curve divides\ the set of all points
of $%
\mathbb{Z}
^{2}$ which are \textquotedblleft exterior\textquotedblright\ to it into
zero, one or two \textquotedblleft connected components\textquotedblright ,
and that these components are infinite.

As an application, we consider curves which are limits of successive
antiderivatives of a complete folding curve. Any such curve is equal to the
closure of its interior. We show that, except in a special case, its
exterior is the union of zero, one or two connected components. In some
cases, its boundary is a fractal.

Finally we prove that, for each finite subcurve $F$ of a complete folding
curve $C$, there exist everywhere in $C$ some subcurves which are parallel
to $F$.

We provide $%
\mathbb{R}
^{2}$ with the euclidian distance defined as $d((x,y),(x^{\prime },y^{\prime
}))=\sqrt{(x^{\prime }-x)^{2}+(y^{\prime }-y)^{2}}$\ for any $x,y,x^{\prime
},y^{\prime }\in 
\mathbb{R}
$.

\bigskip 
\begin{center}
\includegraphics[scale=0.20]{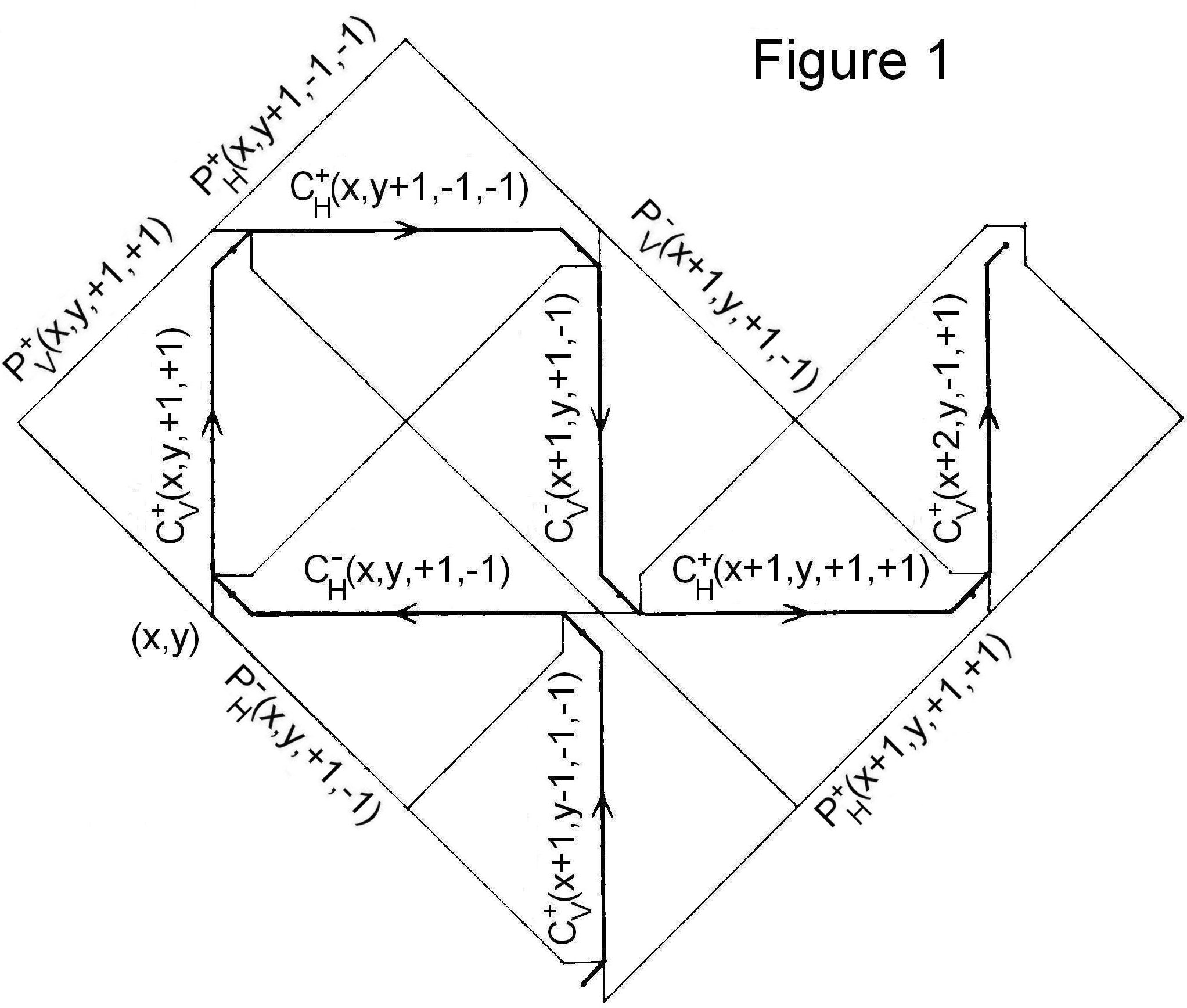}
\end{center}
\bigskip

We fix $\alpha \in 
\mathbb{R}
_{+}^{\ast }$ small compared to $1$. For any $x,y\in 
\mathbb{Z}
$ and any $\zeta ,\eta \in \{+1,-1\}$, we consider (cf. Fig. 1) the \emph{%
segments of curves}

\noindent $C_{H}(x,y,\zeta ,\eta )=[(x+\alpha ,y+\zeta \alpha ),(x+2\alpha
,y)]\cup \lbrack (x+2\alpha ,y),(x+1-2\alpha ,y)]$

$\ \ \ \ \ \ \ \ \ \ \ \ \ \ \ \ \ \cup $ $[(x+1-2\alpha ,y),(x+1-\alpha
,y+\eta \alpha )]$ and

\noindent $C_{V}(x,y,\zeta ,\eta )=[(x+\zeta \alpha ,y+\alpha ),(x,y+2\alpha
)]\cup \lbrack (x,y+2\alpha ),(x,y+1-2\alpha )]$

$\ \ \ \ \ \ \ \ \ \ \ \ \ \ \ \ \ \cup $ $[(x,y+1-2\alpha ),(x+\eta \alpha
,y+1-\alpha )]$.

\noindent We say that $[(x,y),(x+1,y)]$ is the \emph{support} of $%
C_{H}(x,y,\zeta ,\eta )$\ and $[(x,y),(x,y+1)]$ is the \emph{support} of $%
C_{V}(x,y,\zeta ,\eta )$.

We denote by $C_{H}^{+}(x,y,\zeta ,\eta )$ the segment $C_{H}(x,y,\zeta
,\eta )$\ oriented from left to right, and $C_{H}^{-}(x,y,\zeta ,\eta )$ the
segment $C_{H}(x,y,\zeta ,\eta )$\ oriented from right to left. Similarly,
we denote by $C_{V}^{+}(x,y,\zeta ,\eta )$ the segment $C_{V}(x,y,\zeta
,\eta )$\ oriented from bottom to top, and $C_{V}^{-}(x,y,\zeta ,\eta )$ the
segment $C_{V}(x,y,\zeta ,\eta )$\ oriented from top to bottom. From now on,
all the segments considered are oriented.

We associate to $C_{H}^{+}(x,y,\zeta ,\eta )$\ the tile

\noindent $P_{H}^{+}(x,y,\zeta ,\eta )=\{(u,v)\in 
\mathbb{R}
^{2}\mid \left\vert u-(x+1/2)\right\vert +\left\vert v-y\right\vert \leq
1/2\}$

\noindent $\cup \;\{(u,v)\in 
\mathbb{R}
^{2}\mid \sup (\left\vert u-(x+1-\alpha )\right\vert ,\left\vert v-(y+\eta
\alpha )\right\vert )\leq \alpha \}$

\noindent $-\;\{(u,v)\in 
\mathbb{R}
^{2}\mid \sup (\left\vert u-(x+\alpha )\right\vert ,\left\vert v-(y+\zeta
\alpha )\right\vert )<\alpha \}$,

\noindent\ and to $C_{H}^{-}(x,y,\zeta ,\eta )$\ the tile

\noindent $P_{H}^{-}(x,y,\zeta ,\eta )=\{(u,v)\in 
\mathbb{R}
^{2}\mid \left\vert u-(x+1/2)\right\vert +\left\vert v-y\right\vert \leq
1/2\}$

\noindent $\cup \;\{(u,v)\in 
\mathbb{R}
^{2}\mid \sup (\left\vert u-(x+\alpha )\right\vert ,\left\vert v-(y+\zeta
\alpha )\right\vert )\leq \alpha \}$

\noindent $-\;\{(u,v)\in 
\mathbb{R}
^{2}\mid \sup (\left\vert u-(x+1-\alpha )\right\vert ,\left\vert v-(y+\eta
\alpha )\right\vert )<\alpha \}$.

Similarly, we associate to $C_{V}^{+}(x,y,\zeta ,\eta )$\ the tile

\noindent $P_{V}^{+}(x,y,\zeta ,\eta )=\{(u,v)\in 
\mathbb{R}
^{2}\mid \left\vert u-x\right\vert +\left\vert v-(y+1/2)\right\vert \leq
1/2\}$

\noindent $\cup \;\{(u,v)\in 
\mathbb{R}
^{2}\mid \sup (\left\vert u-(x+\eta \alpha )\right\vert ,\left\vert
v-(y+1-\alpha )\right\vert )\leq \alpha \}$

\noindent $-\;\{(u,v)\in 
\mathbb{R}
^{2}\mid \sup (\left\vert u-(x+\zeta \alpha )\right\vert ,\left\vert
v-(y+\alpha )\right\vert )<\alpha \}$,

\noindent\ and to $C_{V}^{-}(x,y,\zeta ,\eta )$\ the tile

\noindent $P_{V}^{-}(x,y,\zeta ,\eta )=\{(u,v)\in 
\mathbb{R}
^{2}\mid \left\vert u-x\right\vert +\left\vert v-(y+1/2)\right\vert \leq
1/2\}$

\noindent $\cup \;\{(u,v)\in 
\mathbb{R}
^{2}\mid \sup (\left\vert u-(x+\zeta \alpha )\right\vert ,\left\vert
v-(y+\alpha )\right\vert )\leq \alpha \}$

\noindent $-\;\{(u,v)\in 
\mathbb{R}
^{2}\mid \sup (\left\vert u-(x+\eta \alpha )\right\vert ,\left\vert
v-(y+1-\alpha )\right\vert )<\alpha \}$.

We say that two segments $C_{1},C_{2}$ are \emph{consecutive} if they just
have one common point and if the end of $C_{1}$ is the beginning of $C_{2}$.
This property is true if and only if the intersection of the associated
tiles consists of one of their four edges (see Fig. 1). The supports of two
consecutive segments form a right angle.

A\emph{\ finite} (resp. \emph{infinite}, \emph{complete}) \emph{curve} is a
sequence $(C_{1},...,C_{n})$ (resp. $(C_{i})_{i\in 
\mathbb{N}
^{\ast }}$, $(C_{i})_{i\in 
\mathbb{Z}
}$) of consecutive segments which all have distinct supports. We identify
two finite curves if they only differ in the beginning of the first segment
and the end of the last one.

The tiles associated to the segments of a curve are nonoverlapping. If we
erase the bumps on their edges, we obtain nonoverlapping square tiles which
cover the same part of $%
\mathbb{R}
^{2}$ if the curve is complete.

We consider that two curves $(C_{1},...,C_{m})$ and $(D_{1},...,D_{n})$ can
be concatenated if their segments have distinct supports and if the end of
the support of $C_{m}$ and the beginning of the support of $D_{1}$ form a
right angle. Then we modify the end of $C_{m}$ and the beginning of $D_{1}$\
in order to make them compatible.

For each finite curve $(C_{i})_{1\leq i\leq n}$ (resp. infinite curve $%
(C_{i})_{i\in 
\mathbb{N}
^{\ast }}$, complete curve $(C_{i})_{i\in 
\mathbb{Z}
}$), we consider the sequence $(\eta _{i})_{1\leq i\leq n-1}$\ (resp. $(\eta
_{i})_{i\in 
\mathbb{N}
^{\ast }}$, $(\eta _{i})_{i\in 
\mathbb{Z}
}$) defined as follows: for each $i$, we write $\eta _{i}=+1$ (resp. $\eta
_{i}=-1$) if we turn left (resp. right) when we pass from $C_{i}$\ to $%
C_{i+1}$. Two curves are associated to the same sequence\ if and only if
they are equivalent modulo a positive isometry.

For each segment of curve $D$, we denote by $\overline{D}$ the segment
obtained from $D$\ by changing the orientation. If a finite curve $%
C=(C_{1},...,C_{n})$ is associated to $S=(\eta _{1},...,\eta _{n-1})$, then $%
\overline{C}=(\overline{C}_{n},...,\overline{C}_{1})$ is associated to $%
\overline{S}=(-\eta _{n-1},...,-\eta _{1})$. If a complete curve $%
C=(C_{i})_{i\in 
\mathbb{Z}
}$ is associated to $S=(\eta _{i})_{i\in 
\mathbb{Z}
}$, then $\overline{C}=(\overline{C}_{-i+1})_{i\in 
\mathbb{Z}
}$ is associated to $\overline{S}=(-\eta _{-i})_{i\in 
\mathbb{Z}
}$.

For each segment of curve $D$ with support $[X,Y]$ and oriented from $X$ to $%
Y$, we call $X,Y$ the \emph{endpoints}, $X$ the \emph{initial point} and $Y$
the \emph{terminal point} of $D$.\ The \emph{initial point} of a curve $%
(C_{1},...,C_{n})$ is the initial point of $C_{1}$, and its \emph{terminal
point} is the terminal point of $C_{n}$. The $\emph{vertices}$ of a curve
are the endpoints of its segments.

We say that two segments of curves, or two curves, are \emph{parallel}
(resp. \emph{opposite}) if they are equivalent modulo a translation (resp. a
rotation of angle $\pi $).

We have $%
\mathbb{Z}
^{2}=M_{1}\cup M_{2}$ and $M_{1}\cap M_{2}=\varnothing $ for $%
M_{1}=\{(x,y)\in 
\mathbb{Z}
^{2}\mid x+y$ odd$\}$ and $M_{2}=\{(x,y)\in 
\mathbb{Z}
^{2}\mid x+y$ even$\}$.\ We denote by $M$ one of these two sets and we
consider, on the one hand the curves with supports of length $1$\ and
vertices in $%
\mathbb{Z}
^{2}$, on the other hand the curves with supports of length $\sqrt{2}$\ and
vertices in $M$.

\bigskip 
\begin{center}
\includegraphics[scale=0.17]{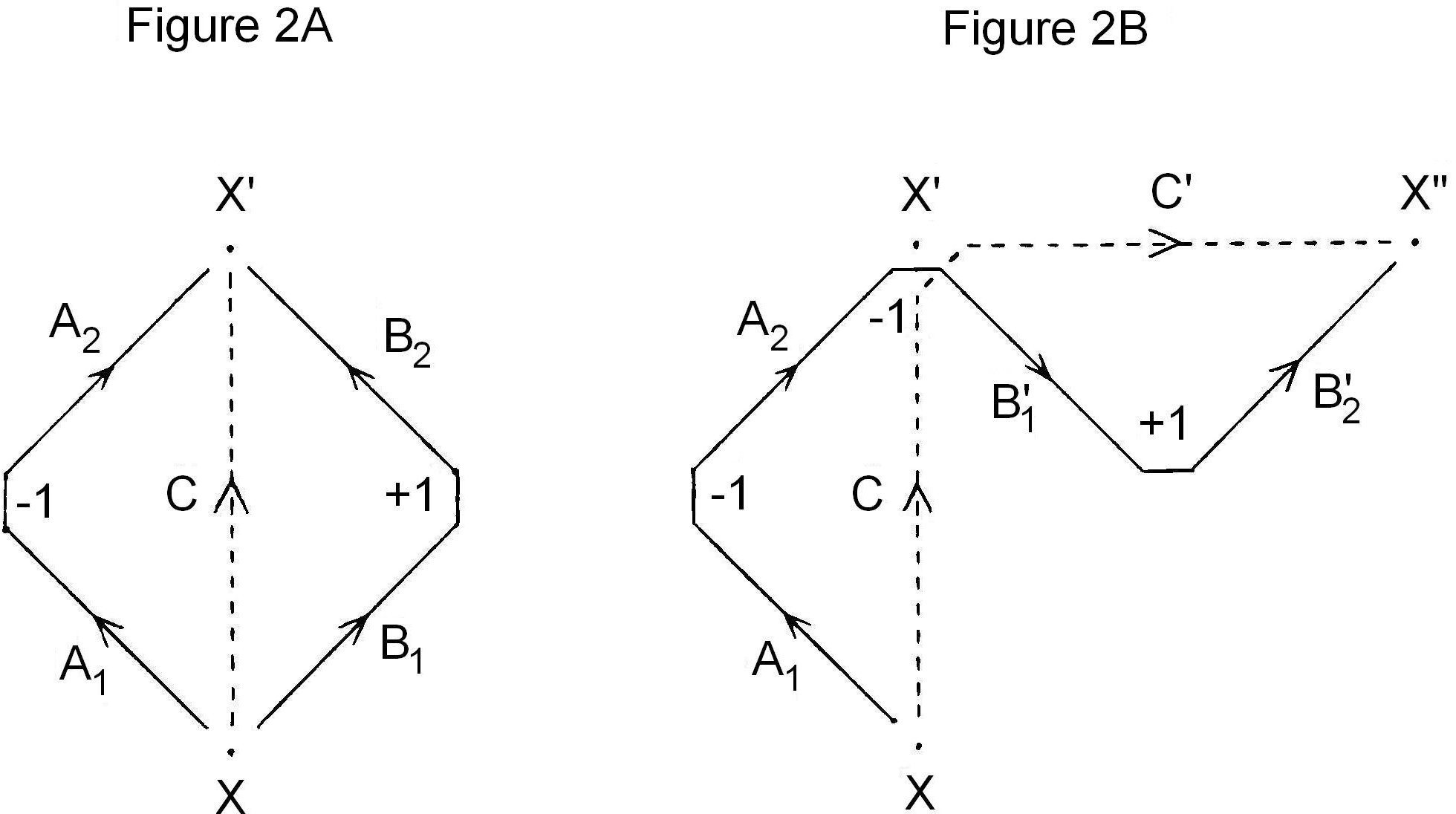}
\end{center}
\bigskip

Let $C$ be a segment of the
second system, let $X$\ be its initial point and let $X^{\prime }$ be its
terminal point. Then, in the first system, there exist two curves $%
(A_{1},A_{2})$ and $(B_{1},B_{2})$, associated to the sequences $(-1)$ and $%
(+1)$, such that $X$ is the initial point of $A_{1}$ and $B_{1}$, and $%
X^{\prime }$ is the terminal point of $A_{2}$ and $B_{2}$ (see Fig. 2A).

Now, consider in the second system a segment $C^{\prime }$ such that $%
(C,C^{\prime })$ is a curve associated to a sequence $(\varepsilon )$\ with $%
\varepsilon \in \{+1,-1\}$. Let $X^{\prime \prime }$ be the terminal point
of $C^{\prime }$. In the first system, denote by $(A_{1}^{\prime
},A_{2}^{\prime })$ and $(B_{1}^{\prime },B_{2}^{\prime })$ the curves
associated to the sequences $(-1)$ and $(+1)$, such that $X^{\prime }$ is
the initial point of $A_{1}^{\prime }$ and $B_{1}^{\prime }$, and $X^{\prime
\prime }$ is the terminal point of $A_{2}^{\prime }$ and $B_{2}^{\prime }$.

Then (see Fig. 2B), $(A_{1},A_{2},B_{1}^{\prime },B_{2}^{\prime })$ and $%
(B_{1},B_{2},A_{1}^{\prime },A_{2}^{\prime })$ are curves associated to $%
(-1,\varepsilon ,+1)$ and $(+1,\varepsilon ,-1)$. Each of these curves has $%
X $, $X^{\prime }$, $X^{\prime \prime }$ among its vertices, and crosses the
curve $(C,C^{\prime })$\ near $X^{\prime }$. Moreover $(A_{1},A_{2},A_{1}^{%
\prime },A_{2}^{\prime })$ and $(B_{1},B_{2},B_{1}^{\prime },B_{2}^{\prime
}) $ are not curves.

For each curve $(C_{1},...,C_{2n})$ (resp. $(C_{i})_{i\in 
\mathbb{N}
^{\ast }}$, $(C_{i})_{i\in 
\mathbb{Z}
}$) of the first system and each curve $(D_{1},...,D_{n})$ (resp. $%
(D_{i})_{i\in 
\mathbb{N}
^{\ast }}$, $(D_{i})_{i\in 
\mathbb{Z}
}$) of the second system, we say that $C$ is an \emph{antiderivative} of $D$
or that $D$ is the \emph{derivative} of $C$ if, for each integer $i$:

\noindent a) if $D_{i+1}$ exists, then $C_{2i+1}$ and $D_{i+1}$ have the
same initial point;

\noindent b) if $D_{i}$ exists, then $C_{2i}$ and $D_{i}$ have the same
terminal point;

\noindent c) if $D_{i}$ and $D_{i+1}$ exist, then $C$ crosses $D$ near the
terminal point of $D_{i}$ (we say that $C$ \emph{alternates} around $D$).

Each curve of the second system has exactly two antiderivatives in the first
one. Each curve $C$ of the first system has at most one derivative in the
second one. If that derivative exists, then the sequence $(\eta
_{1},...,\eta _{2n-1})$ (resp. $(\eta _{i})_{i\in 
\mathbb{N}
^{\ast }}$, $(\eta _{i})_{i\in 
\mathbb{Z}
}$) of elements of $\{-1,+1\}$\ associated to $C$ satisfies $\eta
_{2i+1}=(-1)^{i}\eta _{1}$\ for each integer $i$ such that $\eta _{2i+1}$\
exists. Conversely, if this condition is satisfied, then the \emph{derivative%
} of $C$ is defined by taking for $M$ the set\ $M_{1}$ or $M_{2}$\ which
contains the initial point of $C_{1}$, and replacing each pair of segments $%
(C_{2i-1},C_{2i})$\ with a segment $D_{i}$.

For the definition of the derivative of a complete curve $(C_{i})_{i\in 
\mathbb{Z}
}$, we permit ourselves to change the initial point of indexation, i.e.\ to
replace $(C_{i})_{i\in 
\mathbb{Z}
}$ with $(C_{i+k})_{i\in 
\mathbb{Z}
}$ for an integer $k$. With this convention, the derivative exists if and
only if $\eta _{2i}=(-1)^{i}\eta _{0}$\ for each $i\in 
\mathbb{Z}
$ or $\eta _{2i+1}=(-1)^{i}\eta _{1}$\ for each $i\in 
\mathbb{Z}
$.\ If these two conditions are simultaneously satisfied, we obtain two
different derivatives; in that case, we consider that the derivative is not
defined. This situation never appears for complete folding curves, which
will be considered here.

We define by induction the\emph{\ }$n$\emph{-th} \emph{derivative} $C^{(n)}$%
\ of a curve $C$, with $C^{(0)}=C$ and $C^{(n+1)}$\ derivative of $C^{(n)}$\
for each $n\in 
\mathbb{N}
$, as well as the\emph{\ }$n$\emph{-th} \emph{antiderivatives}. It is
convenient to represent the successive derivatives of a curve $C$ on the
same figure in such a way that $C^{(n)}$ alternates around $C^{(n+1)}$ for
each $n\in 
\mathbb{N}
$\ such that $C^{(n+1)}$ exists. This convention will be used later in the
paper.

For each $n\in 
\mathbb{N}
$, an $n$\emph{-folding curve} is a curve associated to an $n$-folding
sequence. For each $n$-folding sequence obtained by folding $n$ times a
strip of paper, we obtain the associated $n$-folding curve by keeping the
strip folded according to right angles instead of unfolding it completely.

\bigskip 
\begin{center}
\includegraphics[scale=0.18]{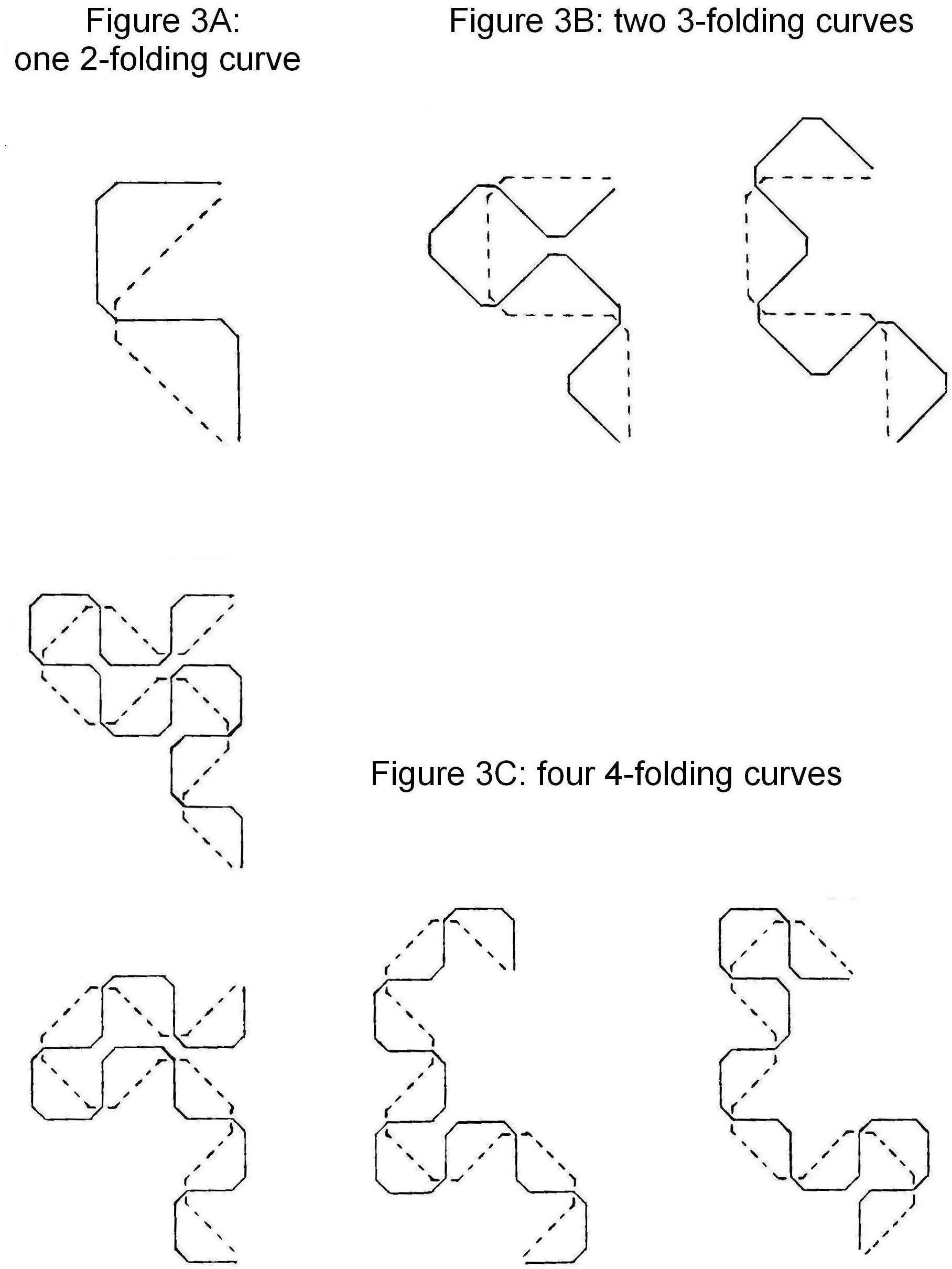}
\end{center}
\bigskip

We see by induction on $n$ that
the $n$-folding curves are the $n$-th antiderivatives of the curves which
consist of one segment. Consequently, up to isometry and up to the
orientation, there exist one $2$-folding curve (cf. Fig. 3A), two $3$%
-folding curves (cf. Fig. 3B), and four $4$-folding curves (cf. Fig. 3C).

We call $\infty $\emph{-folding curve} (resp.\emph{\ finite} \emph{folding
curve}, \emph{complete folding curve}) each curve associated to an $\infty $%
-folding sequence (resp. a finite folding sequence, a complete folding
sequence). Any curve $(C_{i})_{i\in 
\mathbb{N}
^{\ast }}$ (resp. $(C_{i})_{i\in 
\mathbb{Z}
}$) is an $\infty $-folding curve (resp. a complete folding curve) if and
only if it is indefinitely derivable.

The successive antiderivatives of a paperfolding curve, as well as its
successive derivatives if they exist, are also paperfolding curves.

We say that a curve $(C_{1},...,C_{n})$ (resp. $(C_{i})_{i\in 
\mathbb{N}
^{\ast }}$, $(C_{i})_{i\in 
\mathbb{Z}
}$) is \emph{self-avoiding} if we have $C_{i}\cap C_{j}=\emptyset $ for $%
\left\vert j-i\right\vert \geq 2$. Such a curve defines an injective
continuous function from a closed connected subset of $%
\mathbb{R}
$ to $%
\mathbb{R}
^{2}$.\bigskip

\noindent \textbf{Proposition 2.1.} Antiderivatives of self-avoiding curves
are self-avoiding.\bigskip

\noindent \textbf{Proof.} Consider a curve $C$ whose derivative $D$ is
self-avoiding. If $C$ is not self-avoiding, then there exist two segments of 
$C$ which have the same support. These two segments are necessarily coming
from segments of $D$ which have a common endpoint.

In order to prove that this situation is impossible, we consider the
function $\tau $ which is defined on the set of all supports of segments of $%
D$ with $\tau ([(u,v),(u+1,v+\varepsilon )])=+1$ (resp. $-1$) for each $%
(u,v)\in 
\mathbb{Z}
^{2}$ and each $\varepsilon \in \{-1,+1\}$ such that $C$ is above (resp.
below) $D$\ on $[(u,v),(u+1,v+\varepsilon )]$. It suffices to observe that
the equality $\tau ([(u^{\prime },v^{\prime }),(u^{\prime }+1,v^{\prime
}+\varepsilon ^{\prime })])=(-1)^{u^{\prime }-u}\tau
([(u,v),(u+1,v+\varepsilon )])$ is true wherever $\tau $ is defined. In
fact, it is true for the supports of consecutive segments of $D$ because $C$
alternates around $D$, and it is proved in the general case by induction on
the number of consecutive segments between the two segments considered.~~$%
\blacksquare $\bigskip

\noindent \textbf{Corollary 2.2.} Paperfolding curves are
self-avoiding.\bigskip

\noindent \textbf{Proof.} For each integer $n$, each $n$-folding curve is
self-avoiding because it is the $n$-th antiderivative of a self-avoiding
curve which consists of one segment. Each finite folding curve is
self-avoiding since it is a subcurve of an $n$-folding curve for an integer $%
n$. Complete folding curves and $\infty $-folding curves are self-avoiding
because their finite subcurves are self-avoiding.

Another proof is given by [4, Observation 1.11, p. 134].~~$\blacksquare $%
\bigskip

For each self-avoiding curve $C$ and any $x,y\in 
\mathbb{Z}
$, we write:

\noindent $\rho _{C}([(x,y),(x+1,y)]=+1$ (resp. $-1$) if $C$ contains a
segment with the initial point $(x,y)$\ (resp. $(x+1,y)$) and the terminal
point $(x+1,y)$\ (resp. $(x,y)$);

\noindent $\rho _{C}([(x,y),(x,y+1)])=+1$ (resp. $-1$) if $C$ contains a
segment with the initial point $(x,y)$\ (resp. $(x,y+1)$) and the terminal
point $(x,y+1)$\ (resp. $(x,y)$).

There exists $\varepsilon \in \{-1,+1\}$ such that $\rho
_{C}([(x,y),(x+1,y)])=(-1)^{y-x+\varepsilon }$ and $\rho
_{C}([(x,y),(x,y+1)])=(-1)^{y-x+\varepsilon +1}$ wherever $\rho _{C}$\ is
defined. In fact, $\varepsilon $ is the same for the supports of two
consecutive segments, and we see that it is the same for the supports of any
two segments by induction on the number of consecutive segments between
them. We extend the definition of $\rho _{C}$, according to this property,
to the set of all intervals $[(x,y),(x+1,y)]$\ or\ $[(x,y),(x,y+1)]$\ with $%
x,y\in 
\mathbb{Z}
$.

For each self-avoiding curve $C=(C_{i})_{i\in 
\mathbb{Z}
}$ and any pairs $(C_{i},C_{i+1})$, $(C_{j},C_{j+1})$ of consecutive
segments, if $C_{i}$ and $C_{j}$ have the same terminal point, then, by the
property of $\rho _{C}$\ stated above, we turn left when passing from $C_{i}$
to $C_{i+1}$ if and only if we turn left when passing from $C_{j}$ to $%
C_{j+1}$. For each $X\in 
\mathbb{Z}
^{2}$, we write $\sigma _{C}(X)=+1$\ (resp. $-1$) if $X$ is the common
endpoint of two consecutive segments $C_{i}$, $C_{i+1}$ of $C$ and if we
turn left (resp. right) when passing from $C_{i}$ to $C_{i+1}$.

It follows from the definition of derivatives that, for each $k\in 
\mathbb{N}
$:

\noindent a) if $C^{(2k)}$ exists, then there exists $X_{2k}\in 
\mathbb{Z}
^{2}$\ such that the set of all vertices of $C^{(2k)}$ is contained in $%
E_{2k}(C)=X_{2k}+2^{k}%
\mathbb{Z}
^{2}$;

\noindent b) if $C^{(2k+1)}$ exists, then there exists $X_{2k+1}\in 
\mathbb{Z}
^{2}$\ such that the set of all vertices of $C^{(2k+1)}$ is contained in $%
E_{2k+1}(C)=X_{2k+1}+(2^{k},2^{k})%
\mathbb{Z}
+(2^{k},-2^{k})%
\mathbb{Z}
$.

The set $E_{n}(C)$ is defined for each $n\in 
\mathbb{N}
$\ such that $C^{(n)}$ exists.\ If $C^{(n+1)}$ exists, then we have $%
E_{n+1}(C)\subset E_{n}(C)$; we write $F_{n}(C)=E_{n}(C)-E_{n+1}(C)$.

If $S=(\eta _{i})_{i\in 
\mathbb{Z}
}$ is the sequence associated to a complete folding curve $C=(C_{i})_{i\in 
\mathbb{Z}
}$, then, for each $i\in 
\mathbb{Z}
$ and each $n\in 
\mathbb{N}
$, the terminal point of $C_{i}$\ belongs to $E_{n}(C)$\ if and only if $i$\
belongs to $E_{n}(S)$.

The following lemma applies, in particular, to complete folding
curves:\bigskip

\noindent \textbf{Lemma 2.3.} Let $C$ be a derivable self-avoiding complete
curve. Consider a square $Q=[x,x+1]\times \lbrack y,y+1]$ with $x,y\in 
\mathbb{Z}
$. If four vertices of $Q$ are endpoints of segments of $C$, then at least
three segments of $C$, with two of them consecutive, have supports which are
edges of $Q$. If three vertices of $Q$ are endpoints, then the two edges
determined by these vertices are supports of segments of $C$, or neither of
them is a support. If $C$ is derivable twice and if two vertices of $Q$ are
endpoints, then they are necessarily adjacent.\bigskip

\bigskip 
\begin{center}
\includegraphics[scale=0.17]{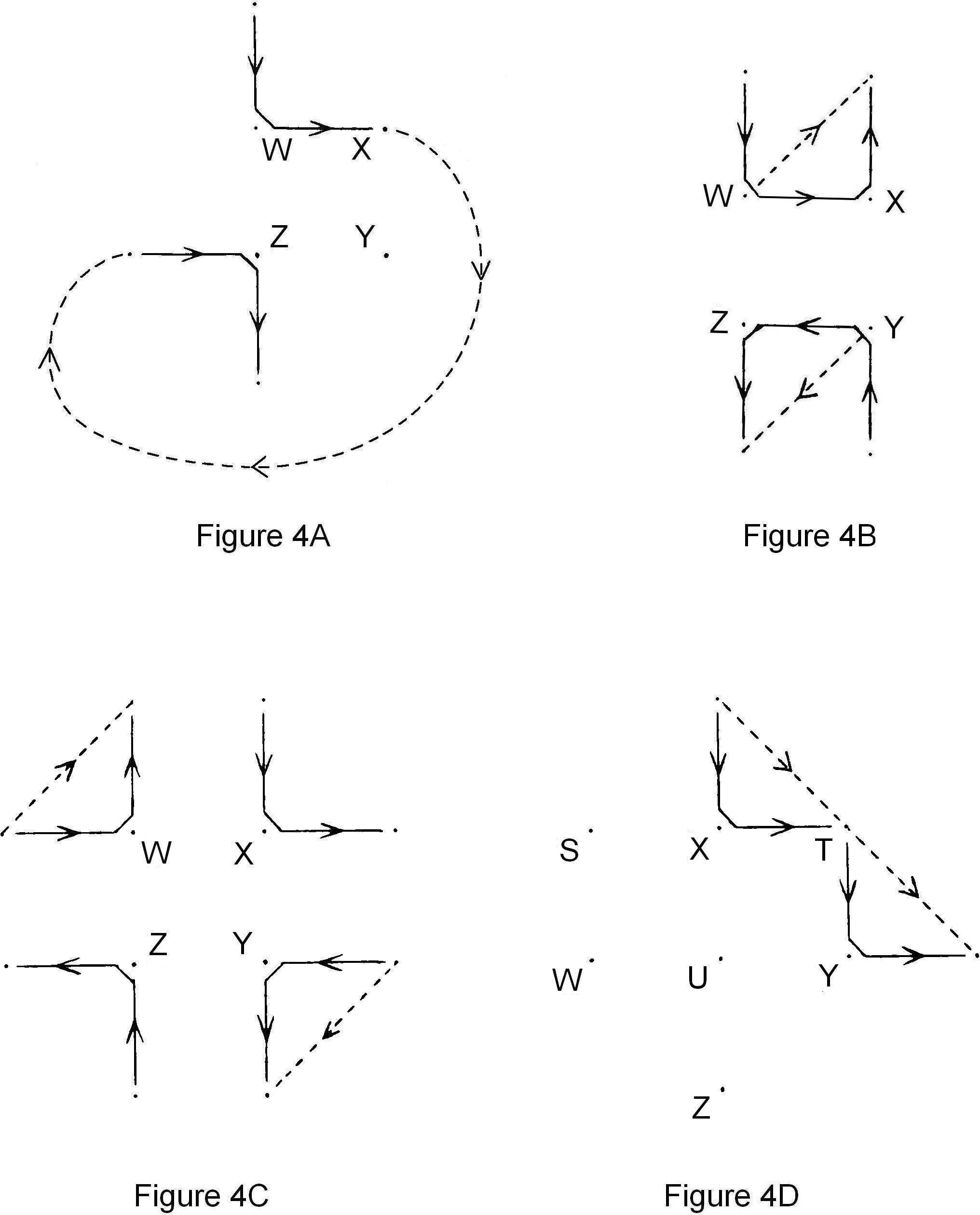}
\end{center}
\bigskip

\noindent \textbf{Proof.} We
denote by $W,X,Y,Z$ the vertices of $Q$ taken consecutively, and we show
that the cases excluded by the Lemma are impossible.

First suppose\ that an edge of $Q$, for instance $WX$, is the support of a
segment of $C$, that a vertex of $Q$ which does not belong to this edge, for
instance $Z$, is a vertex of $C$, and that the edges of $Q$ which contain
this vertex are not supports of segments of $C$. Consider the two pairs of
consecutive segments of $C$ which respectively have $W$ and $Z$\ as a common
endpoint. Then the property of $\rho _{C}$ implies that these two pairs both
have the orientation shown by Figure 4A, or both have the contrary
orientation, which contradicts the connectedness of $C$ (see Fig. 4A).

Then suppose that two opposite edges of $Q$, for instance $WX$ and $YZ$, are
supports of segments of $C$, and that the two preceding segments, as well as
the two following segments, have supports which are not edges of $Q$. Then
the property of $\rho _{C}$ implies that the two sequences of three
consecutive segments formed from these six segments both have the
orientation shown by Figure 4B, or both have the contrary orientation.
Suppose for instance that $X$ and $Z$ belong to $F_{0}(C)$. Then the two
pairs of segments of $C$, extracted from the two sequences, which
respectively have $X$ and $Z$ as a common endpoint, give opposite segments
of $D$, which contradicts the property of $\rho _{D}$ (see Fig. 4B).

Now suppose that $W,X,Y,Z$ are vertices of $C$ and that no edge of $Q$ is
the support of a segment of $C$. Then the property of $\rho _{C}$ implies
that the pairs of segments of $C$ which respectively have $W,X,Y,Z$ as a
common endpoint all have the orientation shown by Figure 4C, or all have the
contrary orientation. Suppose for instance that $W$ and $Y$ belong to $%
F_{0}(C)$. Then the pairs of segments of $C$ which respectively have $W$ and 
$Y$ as a common endpoint give opposite segments of $D$, which contradicts
the property of $\rho _{D}$ (see Fig. 4C).

Finally suppose that $C$ is derivable twice and that only two opposite
vertices of $Q$, for instance $W$ and $Y$, are vertices of $C$. If $W$ and $%
Y $ belong to $F_{0}(C)$, then, as in the previous case, the pairs of
segments of $C$ which respectively have $W$ and $Y$ as a common endpoint
give opposite segments of $D$, which contradicts the property of $\rho _{D}$.

If $W$ and $Y$ belong to $E_{1}(C)$, we consider the two squares of width $%
\sqrt{2}$ which have $WY$ as their common edge. As $W$ and $Y$ are vertices
of $D$, one of the squares has an edge adjacent to $WY$ which is the support
of a segment of $D$. On the other hand, as the center $X$ or $Z$ of that
square is not a vertex of $C$, the edge $WY$ and the opposite edge are not
supports of segments of $D$, which contradicts the two first statements of
the Lemma applied to $D$.~~$\blacksquare $\bigskip

\noindent \textbf{Definitions.} For any $X,Y\in 
\mathbb{Z}
^{2}$, a \emph{path} from $X$ to $Y$ is a sequence ($X_{0},...,X_{n})\subset 
\mathbb{Z}
^{2}$\ with $n\in 
\mathbb{N}
$, $X_{0}=X$, $X_{n}=Y$ and $d(X_{i-1},X_{i})=1$ for $1\leq i\leq n$. For
each complete curve $C$, we call \emph{exterior} of $C$ and we denote by $%
\mathrm{Ext}(C)$ the set of all points of $%
\mathbb{Z}
^{2}$ which are not vertices of $C$. A \emph{connected component} of $%
\mathrm{Ext}(C)$ is a subset $K$ which is maximal for the following
property: any two points of $K$ are connected by a path which only contains
points of $K$.\bigskip

\noindent \textbf{Theorem 2.4.} The exterior $\mathrm{Ext}(C)$ of a complete
folding curve $C$ is the union of $0$, $1$ or $2$ infinite connected
components, and each of these components is the intersection of $\mathrm{Ext}%
(C)$ with one of the $2$ connected components of $%
\mathbb{R}
^{2}-C$.\bigskip

\noindent \textbf{Lemma 2.4.1.} The connected components of $\mathrm{Ext}(C)$%
\ are infinite.\bigskip

\noindent \textbf{Proof of the Lemma.} For each $n\in 
\mathbb{N}
$,\ we have $\mathrm{Ext}(C^{(n)})=\mathrm{Ext}(C)\cap E_{n}(C)$ since each
point of $E_{n}(C)$ is a vertex of $C^{(n)}$ if and only if it is a vertex
of $C$.

If $K$ is a connected component of $\mathrm{Ext}(C)$, then $K\cap E_{n}(C)$
is a union of connected components of $\mathrm{Ext}(C^{(n)})$ for each $n\in 
\mathbb{N}
$. Otherwise, the smallest integer $n$ such that this property is false
satisfies $n\geq 1$, and $E_{n-1}(C)$ contains the consecutive vertices $%
W,X,Y,Z$ of a square of width $(\sqrt{2})^{n-1}$ with $W,Y\in \mathrm{Ext}%
(C^{(n)})$, $W\in K$, $Y\notin K$ and $X,Z$\ vertices of $C^{(n-1)}$, which
contradicts Lemma 2.3 applied to $C^{(n-1)}$.

For each connected component $K$ of $\mathrm{Ext}(C)$ and each $n\in 
\mathbb{N}
$,\ we have $\varnothing \subsetneq K\cap E_{n+1}(C)\subsetneq K\cap
E_{n}(C) $, or $K\cap E_{n}(C)$\ is a union of connected components which
consist of one point; in fact, for any $X,Y\in E_{n}(C)$\ with $d(x,y)=(%
\sqrt{2})^{n}$, we have $X\in F_{n}(C)$\ and $Y\in E_{n+1}(C)$, or $Y\in
F_{n}(C)$\ and $X\in E_{n+1}(C)$. Consequently, in order to prove that $%
\mathrm{Ext}(C)$\ has no finite connected component, it suffices to show
that each $\mathrm{Ext}(C^{(n)})$\ has no connected component which consists
of one point.

Suppose that there exist a twice derivable complete curve $D$ and a point $%
U=(u,v)\in 
\mathbb{Z}
^{2}$ such that $\left\{ U\right\} $ is a connected component of $\mathrm{Ext%
}(D)$. Write $W=(u-1,v)$, $X=(u,v+1)$, $Y=(u+1,v)$ and $Z=(u,v-1)$.

If $U$ belongs to $F_{0}(D)$, then $W,X,Y,Z$\ belong to $E_{1}(D)$ and they
are vertices of $D^{(1)}$. By Lemma 2.3, two consecutive segments of $%
D^{(1)} $ have supports which are edges of $WXYZ$. As $D$\ alternates around 
$D^{(1)} $, it follows that $U$ is a vertex of $D$, whence a contradiction.

If $U$ belongs to $E_{1}(D)$, consider the point $S$ (resp. $T$) which forms
a square with $X,U,W$ (resp. $X,U,Y$). Then $SX$ or $XT$ is the support of a
segment of $D$ since $X$ is a vertex of $D$.

Suppose for instance that $XT$ is the support of a segment of $D$. As $Y$ is
a vertex of $D$ contrary to $U$, Lemma 2.3 implies that $TY$ is also the
support of a segment of $D$. As $X$ and $Y$ belong to $F_{0}(D)$, it follows
from the property of $\rho _{D}$ (see Fig. 4D) that there exist two parallel
segments of $D^{(1)}$ such that $T$ is the terminal point of one of them and
the initial point of the other one, whence a contradiction.~~$\blacksquare $%
\bigskip

\noindent \textbf{Proof of the Theorem.} Write $C=(C_{i})_{i\in 
\mathbb{Z}
}$.\ Consider a connected component $K$ of $\mathrm{Ext}(C)$ and write $%
M=\{X\in 
\mathbb{Z}
^{2}-K\mid d(X,K)=1\}$.

Denote by $\Omega $\ the set of all squares $S=[x,x+1]\times \lbrack y,y+1]$
with $x,y\in 
\mathbb{Z}
$ such that $K$\ contains one or two vertices of $S$. For each $S\in \Omega $%
, if $K$ contains one vertex $X$ of $S$, consider the segment $E_{S}$ of
length $\sqrt{2}$ joining the vertices of $S$ adjacent to $X$. If $K$
contains two vertices of $S$, consider the segment $E_{S}$ of length $1$
joining the two other vertices, which are adjacent by Lemma 2.3.

The endpoints of the segments $E_{S}$ for $S\in \Omega $ belong to $M$. We
are going to prove that each $U\in M$ is an endpoint of exactly two such
segments.\ We consider the points $V,W,X,Y\in 
\mathbb{Z}
^{2}$\ with $d(U,V)=d(U,W)=d(U,X)=d(U,Y)=1$ such that $VWXY$\ is a square, $%
V,W$\ are vertices of $C$, and $X\in K$. We denote by $P,Q,R,S$\ the squares
determined by the pairs of edges $(UV,UW)$, $(UW,UX)$, $(UX,UY)$, $(UY,UV)$.

As $C$\ is connected, the fourth vertex of $P$\ does not belong to $K$, and $%
P$ does not belong to $\Omega $. On the other hand, $Q$\ belongs to $\Omega $%
\ since $X$\ belongs to $K$\ contrary to $U$\ and $W$.

If $Y$ is a vertex of $C$, then $R$\ belongs to $\Omega $\ since $X$\
belongs to $K$\ and $U,Y$\ do not belong to $K$. Moreover, Lemma 2.3 implies
that the fourth vertex of $S$\ is a vertex of $C$, since $UV$\ is the
support of a segment of $C$\ contrary to $UY$. Consequently, $S$\ does not
belong to $\Omega $.

If $Y$ is not a vertex of $C$, then, by Lemma 2.3, the fourth vertex of $R$\
is not a vertex of $C$, and therefore belongs to $K$. Consequently, $Y$\
belongs to $K$, $S$\ belongs to $\Omega $\ and $R$\ does not belong to $%
\Omega $.

Moreover, $U$\ is the common endpoint of $E_{Q}$ and $E_{R}$ if $Y$\ is a
vertex of $C$, and the common endpoint of $E_{Q}$ and $E_{S}$ if $Y$\ is not
a vertex of $C$.

As $K$ is infinite by Lemma 2.4.1, it follows that the segments $E_{S}$ for $%
S\in \Omega $ form an unbounded self-avoiding curve $E$. The vertices of $E$%
\ are the points of $M$. One connected component of $%
\mathbb{R}
^{2}-E$ contains $K$, but contains no point of $C$\ and no point of\ $%
\mathbb{Z}
^{2}$\ which does not belong to $K$.

The points of $M$ taken along $E$ form a sequence $(X_{i})_{i\in 
\mathbb{Z}
}$. For each $X\in M$, denote by $r(X)$ the unique integer $j$ such that $X$
is the terminal point of $C_{j}$ and the initial point of $C_{j+1}$. Suppose
for instance $r(X_{0})<r(X_{1})$. Then, using the connectedness of $C$, we
see by induction on $i$ that $r(X_{i})<r(X_{i+1})$ for each\ $i\in 
\mathbb{Z}
$.

Now suppose that there exists a connected component $L\neq K$ of $\mathrm{Ext%
}(C)$ such that $K$ and $L$ are contained in the same connected component of 
$%
\mathbb{R}
^{2}-C$. Then there exists\ $i\in 
\mathbb{Z}
$ such that $L$ is contained in the loop formed by $C$ between $X_{i}$ and $%
X_{i+1}$, and $L$ is finite contrary to Lemma 2.4.1.~~$\blacksquare $\bigskip

Now, we apply Theorem 2.4 to limits of complete folding curves. We give less
details in this part, which will not be used in the remainder of the paper.

We consider some complete folding curves $C_{n}=(C_{n,p})_{p\in 
\mathbb{Z}
}$\ with $C_{n}=C_{n+1}^{(1)}$\ for each $n\in 
\mathbb{N}
$.\ We suppose that the curves $C_{n}$\ are represented on the same figure
in such a way that:

\noindent - $C_{0}$ has vertices in $%
\mathbb{Z}
^{2}$\ and supports of length $1$;

\noindent - all the segments $C_{n,1}$\ have the same initial point;

\noindent - $C_{n+1}$\ alternates around $C_{n}$\ for each $n\in 
\mathbb{N}
$.

\noindent We denote by $L$ the limit of the curves $C_{n}$ considered as
representations of functions from $%
\mathbb{R}
$\ to $%
\mathbb{R}
^{2}$.

The curve $L$ is associated to a function from $%
\mathbb{R}
$\ to $%
\mathbb{R}
^{2}$\ which is continuous everywhere, but derivable nowhere. Moreover, $L$\
is closed as a subset of $%
\mathbb{R}
^{2}$. By Theorem 3.1 below, $L$\ contains arbitrarily large open balls. It
follows from Proposition 2.6 and Theorem 3.1 that $L$ is equal to the
closure of its interior.

Now, we consider the complete folding sequences $(\eta _{n,p})_{p\in 
\mathbb{Z}
}$\ associated to the curves $C_{n}$. We say that $(\eta _{n,1})_{n\in 
\mathbb{N}
}$\ is \emph{ultimately alternating} if we have $\eta _{n,1}=(-1)^{n}$\ for $%
n$\ large enough, or\ $\eta _{n,1}=(-1)^{n+1}$\ for $n$\ large
enough.\bigskip

\noindent \textbf{Corollary 2.5.} If $(\eta _{n,1})_{n\in 
\mathbb{N}
}$\ is not ultimately alternating, then $%
\mathbb{R}
^{2}-L$ is the union of $0$, $1$\ or $2$\ unbounded connected
components.\bigskip

\noindent \textbf{Proof.} For each $n\in 
\mathbb{N}
$, we have $\mathrm{Ext}(C_{n})\subset 
\mathbb{R}
^{2}-L$ since $(\eta _{m,1})_{m\in 
\mathbb{N}
}$\ is not ultimately alternating. By Theorem 2.4, it suffices to show
that,\ for each $n\in 
\mathbb{N}
$, any $X,X^{\prime }\in $ $\mathrm{Ext}(C_{n})$ belong to the same
connected component of $%
\mathbb{R}
^{2}-L$ if they belong to the same connected component of $\mathrm{Ext}%
(C_{n})$. Moreover, it is enough to prove this property when $d(X,X^{\prime
})=(1/\sqrt{2})^{n}$.

We consider some distinct points $Y,Y^{\prime },Z,Z^{\prime },U,V$\ with $%
YY^{\prime }$\ and $ZZ^{\prime }$\ parallel to $XX^{\prime }$\ such that $%
XX^{\prime }Y^{\prime }Y$ (resp. $XX^{\prime }Z^{\prime }Z$) is a square of
center $U$ (resp. $V$). Then $XX^{\prime }$, $XY$, $XZ$, $X^{\prime
}Y^{\prime }$, $X^{\prime }Z^{\prime }$ are not supports of segments of $%
C_{n}$. Consequently, at least one of the points $U,V$ is not a vertex of $%
C_{n+1}$.

If $U$ ( resp. $V$) is not a vertex of $C_{n+1}$, then the open triangle $%
XUX^{\prime }$\ (resp. $XVX^{\prime }$) is contained in $%
\mathbb{R}
^{2}-L$. It follows that $X$\ and $X^{\prime }$ belong to the same connected
component of $%
\mathbb{R}
^{2}-L$.~~$\blacksquare $\bigskip

\noindent \textbf{Remark.} Let $L$ be the limit of the curves $C_{n}$, where 
$C_{0}$\ is the curve $C$\ of Example 3.13 below and $\eta _{n,1}=(-1)^{n+1}$%
\ for $n\geq 1$. Then we have $[(0,1),(2,1)]\subset L$ even though $(1,1)$
is not a vertex of $C_{0}$. It follows that $(1,0)$\ belongs to a bounded
connected component of $%
\mathbb{R}
^{2}-L$. It can be proved that $%
\mathbb{R}
^{2}-L$ has infinitely many such components.\bigskip

\noindent \textbf{Examples. }The limit curve $L$\ obtained from the curve $C$
of Example 3.13 by writing $\eta _{n,1}=+1$ for $n\geq 1$ is called a \emph{%
dragon curve}. It follows from [6] that the interior of $L$ is a union of
countably many bounded connected components. According to [2], the boundary
of $L$ is a fractal.

On the other hand, Example 3.8 gives a case with $L=%
\mathbb{R}
^{2}$. Using a similar construction, we can obtain $L$ such that its
boundary is a line. If $L$ is one of the limit curves obtained from the
curves of Example 3.14 by writing $\eta _{n,1}=(-1)^{n+1}$\ for $n\geq 1$,
then its boundary is the union of two intersecting lines, or the union of
two half-lines with a common endpoint.\bigskip

The Proposition below is not a priori obvious since two curves associated to
the same folding sequence are not necessarily parallel.\bigskip

\noindent \textbf{Proposition 2.6.} For each complete folding curve $%
C=(C_{h})_{h\in 
\mathbb{Z}
}$, for each $n\in 
\mathbb{N}
$ and for any $i,j\in 
\mathbb{Z}
$, $(C_{j+1},...,C_{j+88n})$\ contains a curve which is parallel and another
curve which is opposite to $(C_{i+1},...,C_{i+n})$.\bigskip

For the proof of Proposition 2.6, we fix $C,n,i,j$ and we write $%
A=(C_{i+1},...,C_{i+n})$. We consider the integer $m$ such that $%
2^{m-1}<n\leq 2^{m}$, and the sequence $S=(\eta _{h})_{h\in 
\mathbb{Z}
}$ associated to $C$.\bigskip

\noindent \textbf{Lemma 2.6.1.} There exists $k\in E_{m+2}(S)$\ such that $%
(C_{k+1},...,C_{k+2^{m+2}})$ contains a subcurve which is parallel or
opposite to $A$.\bigskip

\noindent \textbf{Proof of the Lemma.} We can suppose that there exists $%
l\in E_{m+2}(S)$ such that $i+1\leq l\leq i+n-1$ since, otherwise, there
exists $k\in E_{m+2}(S)$\ such that $A\subset (C_{k+1},...,C_{k+2^{m+2}})$.

Then, for each $h\in \{i+1,...,i+n-1\}-\{l\}$, we have $h\in 
\mathbb{Z}
-E_{n}$, and therefore $\eta _{h}=\eta _{h+2^{m+1}+r.2^{m+2}}$ for each $%
r\in 
\mathbb{Z}
$. We also have $\eta _{l+2^{m+1}+r.2^{m+2}}=(-1)^{r}\eta _{l+2^{m+1}}$ for
each $r\in 
\mathbb{Z}
$. Consequently, there exists $r\in 
\mathbb{Z}
$ such that $\eta _{l}=\eta _{l+2^{m+1}+r.2^{m+2}}$, and therefore $\eta
_{h}=\eta _{h+2^{m+1}+r.2^{m+2}}$\ for $i+1\leq h\leq i+n-1$.

Now the supports of $C_{i+1}$ and $C_{i+1+2^{m+1}+r.2^{m+2}}$\ are parallel
or opposite since $(i+1+2^{m+1}+r.2^{m+2})-(i+1)$\ is even. It follows that

\noindent $(C_{i+1+2^{m+1}+r.2^{m+2}},...,C_{i+n+2^{m+1}+r.2^{m+2}})$ is
parallel or opposite to $A$.~~$\blacksquare $\bigskip

\noindent \textbf{Proof of the Proposition.} By Lemma 2.6.1, we can suppose
that there exists $k\in E_{m+2}(S)$ such that $k\leq i$ and $k+2^{m+2}\geq
i+n$. For each $r\in 
\mathbb{Z}
$, we write $C_{r}^{\ast }=(C_{k+1+r.2^{m+2}},...,C_{k+(r+1).2^{m+2}})$. We
consider the $(m+2)$-th derivative $D=(D_{r})_{r\in 
\mathbb{Z}
}$\ of $C$, indexed in such a way that, for each $r\in 
\mathbb{Z}
$, the initial point of $D_{r}$ is the initial point of $C_{k+1+r.2^{m+2}}$.

For any $r,s\in 
\mathbb{Z}
$, we have

\noindent $(\eta _{k+1+r.2^{m+3}},...,\eta _{k+2^{m+2}-1+r.2^{m+3}})=(\eta
_{k+1+s.2^{m+3}},...,\eta _{k+2^{m+2}-1+s.2^{m+3}})$.

\noindent Consequently, $C_{2r}^{\ast }$ and $C_{2s}^{\ast }$\ are parallel
(resp. opposite) if and only if $D_{2r}$\ and $D_{2s}$\ are parallel (resp.
opposite). The same property is true for the copies of $A$\ contained in $%
C_{2r}^{\ast }$ and $C_{2s}^{\ast }$.

As we have $88n>11.2^{m+2}$, there exists $r\in 
\mathbb{Z}
$ such that $(C_{j+1},...,C_{j+88n})$\ contains $(C_{2r}^{\ast
},...,C_{2r+8}^{\ast })$. Moreover, $(D_{2r},...,D_{2r+8})$ is necessarily
contained in a $4$-folding curve. It follows (see Fig. 3C) that there exist $%
p,q\in \{0,1,2,3,4\}$\ such that $D_{2r+p}$ and $D_{2r+q}$ are opposite.
Then $C_{2r+p}^{\ast }$ and $C_{2r+q}^{\ast }$ are also opposite, as well as
the copies of $A$ that they contain. Consequently, one of these copies is
parallel to $A$.~~$\blacksquare $\bigskip

\textbf{3. Coverings of $%
\mathbb{R}
^{2}$ by sets of disjoint complete folding curves.}\bigskip

Consider $\Omega =%
\mathbb{R}
^{2}$ or $\Omega =[a,b]\times \lbrack c,d]$\ with\ $a,b,c,d\in 
\mathbb{Z}
$, $b\geq a+1$ and $d\geq c+1$. Let $\mathcal{E}$\ be a set of segments with
supports contained in $\Omega $. We say that $\mathcal{E}$\ is a \emph{%
covering} of $\Omega $\ if it satisfies the following conditions:

\noindent 1) each interval $[(x,y),(x+1,y)]$ or $[(x,y),(x,y+1)]$ contained
in $\Omega $, with $x,y\in 
\mathbb{Z}
$, is the support of a unique segment of $\mathcal{E}$;

\noindent 2) if two distinct non consecutive segments of $\mathcal{E}$ have
a common endpoint $X$, then they can be completed into pairs of consecutive
segments, with all four segments having distinct supports which contain $X$.

\noindent The set $\mathcal{E}$\ is a covering of $%
\mathbb{R}
^{2}$ (resp. $[a,b]\times \lbrack c,d]$) if and only if the tiles associated
to the segments belonging to $\mathcal{E}$\ form a tiling of $%
\mathbb{R}
^{2}$ (resp. a patch which covers $[a,b]\times \lbrack c,d]$).

For each covering $\mathcal{E}$ of $%
\mathbb{R}
^{2}$, if each finite sequence of consecutive segments belonging to $%
\mathcal{E}$\ is a folding curve, then $\mathcal{E}$\ is a \emph{covering of 
}$%
\mathbb{R}
^{2}$\emph{\ by complete folding curves}, in the sense that $\mathcal{E}$ is
a union of disjoint complete folding curves.

We say that a covering of $%
\mathbb{R}
^{2}$ by complete folding curves \emph{satisfies the local isomorphism
property} if the associated tiling satisfies the local isomorphism property.
Two such coverings are said to be \emph{locally isomorphic} if the
associated tilings are locally isomorphic.

It often happens that a covering of $%
\mathbb{R}
^{2}$ by complete folding curves satisfies the local isomorphism property.
In particular, we show that this property is satisfied by any covering of $%
\mathbb{R}
^{2}$ by $1$ complete folding curve, or by $2$ complete folding curves
associated to the \textquotedblleft positive\textquotedblright\ folding
sequence. Two complete folding curves associated to the \textquotedblleft
alternating\textquotedblright\ folding sequence do not give a covering of $%
\mathbb{R}
^{2}$, but we prove that a covering of $%
\mathbb{R}
^{2}$ which satisfies the local isomorphism property is obtained naturally
from these $2$ curves by adding $4$ other complete folding curves.

We characterize the coverings of $%
\mathbb{R}
^{2}$ by sets of complete folding curves which satisfy the local isomorphism
property, and the pairs of locally isomorphic such coverings. We show that
each complete folding curve can be completed in a quasi-unique way into such
a covering and that, for each complete folding sequence $S$, there exists a
covering of $%
\mathbb{R}
^{2}$ by a complete folding curve associated to a sequence which is locally
isomorphic to $S$.

Finally, we prove that each complete folding curve covers a
\textquotedblleft significant\textquotedblright\ part of $%
\mathbb{R}
^{2}$. In that way, we show that the maximum number of disjoint complete
folding curves in $%
\mathbb{R}
^{2}$, and therefore the maximum number of complete folding curves in a
covering of $%
\mathbb{R}
^{2}$, is at most $24$. We also prove that such a covering cannot contain
more than $6$ curves if it satisfies the local isomorphism property.

The following result will be used frequently in the proofs:\bigskip

\noindent \textbf{Theorem 3.1.} There exists a fonction $f:%
\mathbb{N}
\rightarrow 
\mathbb{N}
$\ with exponential growth such that, for each integer $n$, each $n$-folding
curve contains a covering of a square $[x,x+f(n)]\times \lbrack y,y+f(n)]$
with $x,y\in 
\mathbb{Z}
$.\bigskip

\bigskip 
\begin{center}
\includegraphics[scale=0.15]{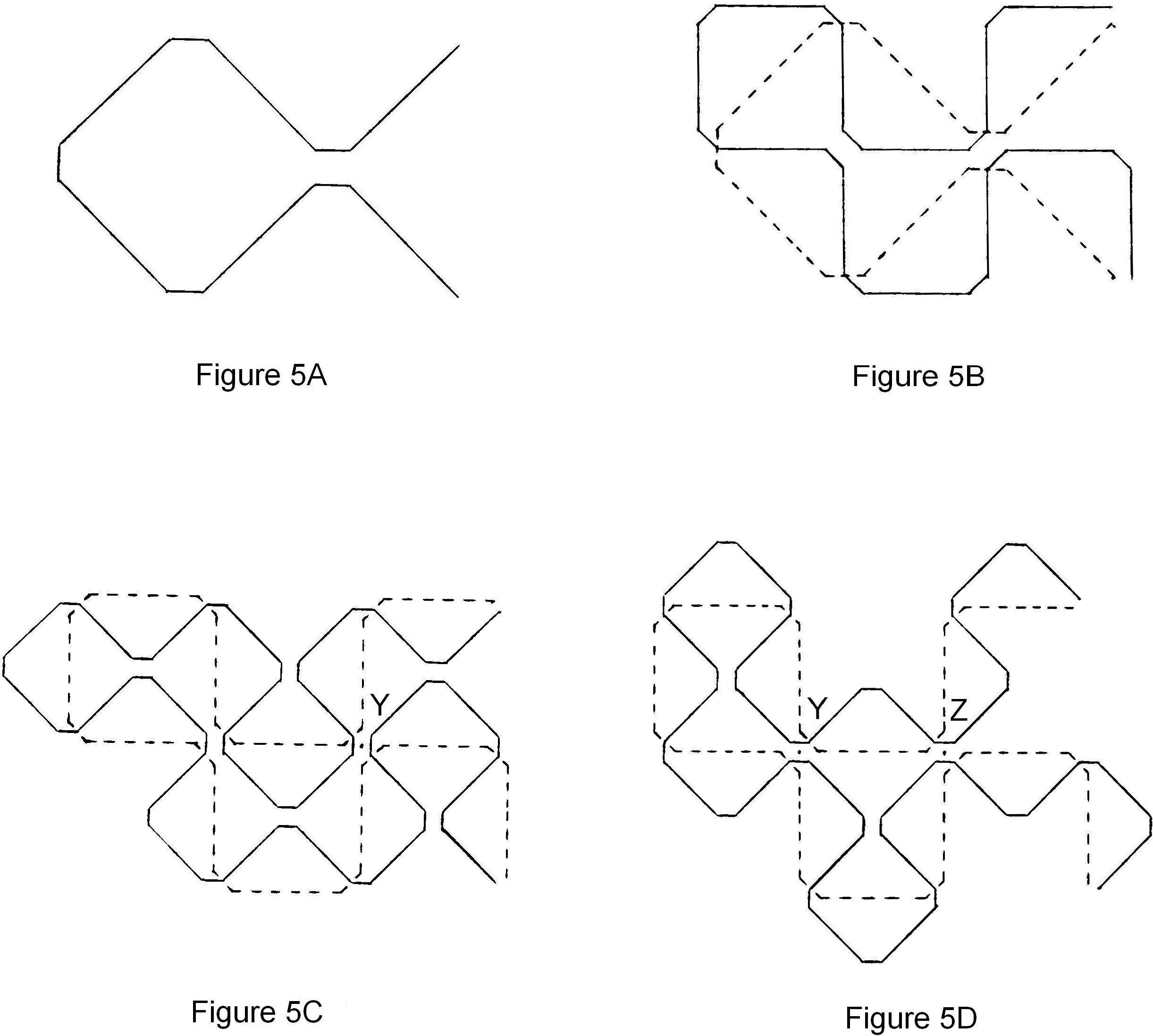}
\end{center}
\bigskip

\noindent \textbf{Proof.} By
Figure 3C, each $4$-folding curve contains a copy, up to isometry and modulo
the orientation, of the curve given by Figure 5A. Consequently, each $5$%
-folding curve contains a copy of the curve given by Figure 5B, and each $6$%
-folding curve contains a copy of one of the two curves given by Figures 5C
and 5D.

For each $X=(x,y)\in 
\mathbb{Z}
^{2}$ and each $k\in 
\mathbb{N}
^{\ast }$, we write $L(X,k)=\{(u,v)\in 
\mathbb{R}
^{2}\mid \left\vert u-x\right\vert +\left\vert v-y\right\vert \leq k\}$. We
say that a folding curve $C$ covers $L(X,k)$ if each interval $%
[(u,v),(u+1,v)]$ or $[(u,v),(u,v+1)]$ with $u,v\in 
\mathbb{Z}
$ contained in $L(X,k)$ is the support of a segment of $C$.

The curve in Figure 5C covers $L(Y,2)$, and each of its two antiderivatives
covers $L(X,2)$ where $X$ is the point corresponding to $Y$. Each
antiderivative of the curve in Figure 5D covers $L(W,2)$ or $L(X,2)$ where $%
W $ and $X$ are the points corresponding to $Y$ and $Z$. Consequently, each $%
7$-folding curve covers an $L(U,2)$.

For each $k\in 
\mathbb{N}
^{\ast }$, if a folding curve $C$ covers an $L(X,k)$, then each second
antiderivative $D$ of $C$ covers an $L(Y,2k-1)$. More precisely, if we put $%
D $ on the figure containing $C$, then $D$ covers $L(X,k-1/2)$, and we
obtain a curve which covers an $L(Y,2k-1)$ when we apply a homothety of
ratio $2$\ in order to give the length $1$ to the supports of the segments
of $D$.

We see by induction on $n$ that each $(2n+7)$-folding curve covers an $%
L(X,2^{n}+1)$.~~$\blacksquare $\bigskip

For each complete folding curve $C=(C_{i})_{i\in 
\mathbb{Z}
}$, there are two possibilities:

\noindent - either all the segments $C_{i}$ for $i\in E_{1}$\ have
horizontal supports, and we say that $C$ has the \emph{type }$H$;

\noindent - or all the segments $C_{i}$ for $i\in E_{1}$\ have vertical
supports, and we say that $C$ has the \emph{type }$V$.\bigskip

\noindent \textbf{Theorem 3.2.} Let $\mathcal{C},\mathcal{D}$ be coverings
of $%
\mathbb{R}
^{2}$ by sets of complete folding curves which satisfy the local isomorphism
property. Then $\mathcal{C}$ and $\mathcal{D}$ are locally isomorphic if and
only if each curve of $\mathcal{C}$ and each curve of $\mathcal{D}$ have the
same type and have locally isomorphic associated sequences.\bigskip

\noindent \textbf{Remark.} In particular, if a covering of $%
\mathbb{R}
^{2}$ by a set of complete folding curves satisfies the local isomorphism
property, then all the curves have the same type and have locally isomorphic
associated sequences.\bigskip

\noindent \textbf{Proof of the Theorem.} First we show that the condition is
necessary. Let $C=(C_{i})_{i\in 
\mathbb{Z}
}$ be a curve of $\mathcal{C}$ and let $D=(D_{i})_{i\in 
\mathbb{Z}
}$ be a curve of $\mathcal{D}$. Consider a finite subcurve $F$ of $C$.

As $\mathcal{C}$ satisfies the local isomorphism property, there exists an
integer $n$ such that each square $[x,x+n]\times \lbrack y,y+n]$ with $%
x,y\in 
\mathbb{Z}
$ contains a subcurve of a curve of $\mathcal{C}$ which is parallel to $F$.
As $\mathcal{C}$ is locally isomorphic to $\mathcal{D}$, each square $%
[x,x+n]\times \lbrack y,y+n]$ with $x,y\in 
\mathbb{Z}
$ also contains a subcurve of a curve of $\mathcal{D}$ which is parallel to $%
F$. According to Theorem 3.1, there exist $x,y\in 
\mathbb{Z}
$ such that $[x,x+n]\times \lbrack y,y+n]$ is covered by $D$. It follows
that $D$ contains a subcurve which is parallel to $F$. In particular, the
folding sequence associated to $F$ is a subword of the folding sequence
associated to $D$.

Now, take for $F$ a $3$-folding subcurve $(C_{i+1},...,C_{i+8})$ of $C$, and
consider $j\in 
\mathbb{Z}
$ such that $(D_{j+1},...,D_{j+8})$ is parallel to $F$. Then we have $i\in
E_{1}(C)$ and $j\in E_{1}(D)$. It follows that $C$ and $D$ have the same
type.

It remains to be proved that the condition is sufficient. We show that, for
each finite set $\mathcal{E}$ of tiles, if $\mathcal{C}$ contains the image
of $\mathcal{E}$\ under a translation, then $\mathcal{D}$ contains the image
of $\mathcal{E}$\ under a translation; the converse can be proved in the
same way.

As $\mathcal{C}$ satisfies the local isomorphism property, there exists an
integer $m$ such that each square $[a,a+m]\times \lbrack b,b+m]$ contains a
set of tiles of $\mathcal{C}$ which is the image of $\mathcal{E}$ under a
translation.\ By Theorem 3.1, there exists an integer $n$ such that each $n$%
-folding curve contains a covering of a square $[a,a+m]\times \lbrack b,b+m]$%
. For such an $n$, each $n$-folding subcurve of a curve of $\mathcal{C}$
contains the image of $\mathcal{E}$ under a translation.

We can take $n\geq 3$. Then each $n$-folding subcurve $F$ of a curve of $%
\mathcal{D}$ is parallel or opposite to a subcurve of a curve of $\mathcal{C}
$ since each curve of $\mathcal{C}$ and each curve of $\mathcal{D}$ have the
same type and have locally isomorphic associated sequences. Now, it follows
from Proposition 2.6 applied to $\mathcal{C}$ that such an $F$ is parallel
to a subcurve of a curve of $\mathcal{C}$, and therefore contains the image
of $\mathcal{E}$ under a translation.~~$\blacksquare $\bigskip

For any disjoint complete folding curves $C,D$, we call \emph{boundary
between }$C$\emph{\ and }$D$ the set of all points of $%
\mathbb{Z}
^{2}$ which are vertices of two segments of $C$ and two segments of $D$%
.\bigskip

\noindent \textbf{Proposition 3.3.} Let $n\geq 1$\ be an integer and let $%
\mathcal{C}$ be a covering of $%
\mathbb{R}
^{2}$ by a set of complete folding curves which satisfies the local
isomorphism property. Then the curves of $\mathcal{C}$ define the same $%
E_{n} $, and their $n$-th derivatives form a covering of $%
\mathbb{R}
^{2}$ by a set of complete folding curves which satisfies the local
isomorphism property. If the boundary between two curves $C,D\in \mathcal{C}$
is nonempty, then the boundary between $C^{(n)}$ and $D^{(n)}$ is
nonempty.\bigskip

\noindent \textbf{Proof.} By induction, it suffices to show the Proposition
for $n=1$. Consequently, it suffices to prove that, for each covering $%
\mathcal{C}$ of $%
\mathbb{R}
^{2}$ by a set of complete folding curves which satisfies the local
isomorphism property, if the boundary between two curves $C,D\in \mathcal{C}$
is nonempty, then we have $E_{1}(C)=E_{1}(D)$, the curves $C^{(1)}$ and $%
D^{(1)}$ are disjoint and the boundary between $C^{(1)}$ and $D^{(1)}$ is
nonempty. Then each point of $E_{1}(\mathcal{C})$ will be an endpoint of $4$
segments of $\mathcal{C}^{(1)}=\{C^{(1)}\mid C\in \mathcal{C}\}$ since it is
an endpoint of $4$ segments of $\mathcal{C}$, and $\mathcal{C}^{(1)}$\ will
satisfy the local isomorphism property like $\mathcal{C}$.

We write $C=(C_{i})_{i\in 
\mathbb{Z}
}$ and $D=(D_{i})_{i\in 
\mathbb{Z}
}$. We denote by $S=(\zeta _{i})_{i\in 
\mathbb{Z}
}$ and $T=(\eta _{i})_{i\in 
\mathbb{Z}
}$ the associated sequences. We consider a point $X$ which belongs to the
boundary between $C$ and $D$. We write $A=(C_{i-4},...,C_{i+3})$ and $%
B=(D_{j-4},...,D_{j+3})$, where $i$ (resp. $j$) is the integer such that $X$
is the common endpoint of $C_{i}$ and $C_{i+1}$\ (resp. $D_{j}$ and $D_{j+1}$%
). As $\mathcal{C}$ satisfies the local isomorphism property, it follows
from Theorem 3.1 applied to $C$ that there exist a translation $\tau $ of $%
\mathbb{R}
^{2}$ and two sequences $A^{\prime }=(C_{k-4},...,C_{k+3})$ and $B^{\prime
}=(C_{l-4},...,C_{l+3})$\ such that $\tau (A)=A^{\prime }$ and $\tau
(B)=B^{\prime }$.

If $X$ belongs to $F_{0}(C)$, then we have $i\in F_{0}(S)$, and therefore $%
\zeta _{i-4}=-\zeta _{i-2}=\zeta _{i}=-\zeta _{i+2}$, which implies $\zeta
_{k-4}=-\zeta _{k-2}=\zeta _{k}=-\zeta _{k+2}$ and $k\in F_{0}(S)$.
Consequently, we have $\tau (X)\in F_{0}(C)$ and $l\in F_{0}(S)$, which
implies $\zeta _{l-4}=-\zeta _{l-2}=\zeta _{l}=-\zeta _{l+2}$ and $\eta
_{j-4}=-\eta _{j-2}=\eta _{j}=-\eta _{j+2}$. It follows that $X$ belongs to $%
F_{0}(D)$. We show in the same way that $X$ belongs to $F_{0}(C)$ if it
belongs to $F_{0}(D)$. Consequently, we have $F_{0}(C)=F_{0}(D)$.

If $X$ belongs to $F_{0}(C)=F_{0}(D)$, then the segment of $C^{(1)}$
obtained from $(C_{i},C_{i+1})$ and the segment of $D^{(1)}$ obtained from $%
(D_{j},D_{j+1})$ have supports which are opposite edges of a square of
center $X$ and width $\sqrt{2}$. Moreover, the segments of $C^{(1)}$
obtained from $(C_{k},C_{k+1})$ and $(C_{l},C_{l+1})$ have supports which
are opposite edges of a square of center $\tau (X)$ and width $\sqrt{2}$. By
Lemma 2.3, a third edge of the second square is the support of a segment of $%
C^{(1)}$ obtained from one of the pairs $(C_{k-2},C_{k-1})$, $%
(C_{k+2},C_{k+3})$, $(C_{l-2},C_{l-1})$, $(C_{l+2},C_{l+3})$. Consequently,
a third edge of the first square is the support of a segment of $C^{(1)}$
obtained from one of the pairs $(C_{i-2},C_{i-1})$, $(C_{i+2},C_{i+3})$, or
the support of a segment of $D^{(1)}$ obtained from one of the pairs $%
(D_{j-2},D_{j-1})$, $(D_{j+2},D_{j+3})$. In both cases, $C^{(1)}$ and $%
D^{(1)}$ have a common vertex.

If $X$ belongs to $E_{1}(C)=E_{1}(D)$, then $X$ is a common vertex of $%
C^{(1)}$ and $D^{(1)}$. Moreover, the two segments of $C^{(1)}$ and the two
segments of $D^{(1)}$ which have $X$ as an endpoint all have different
supports since they are the images under $\tau ^{-1}$ of the four segments
of $C^{(1)}$ which have $\tau (X)$ as an endpoint. As this property is true
for each point of the boundary between $C$\ and $D$\ which belongs to $%
E_{1}(C)=E_{1}(D)$, the curves $C^{(1)}$ and $D^{(1)}$ are disjoint.~~$%
\blacksquare $\bigskip

\noindent \textbf{Remark.} If the boundary between two curves $C,D\in 
\mathcal{C}$ is finite, then it contains a point of $E_{\infty }$.
Otherwise, for $n$ large enough, it would contain no point of $E_{n}$, and
the boundary between $C^{(n)}$ and $D^{(n)}$ would be empty, contrary to
Proposition 3.3.\bigskip

If two disjoint complete folding curves $C,D$ have the same type and if $%
E_{1}(C)=E_{1}(D)$, then we have $\sigma _{C}(X)=\sigma _{D}(X)$\ for each
point $X$ of the boundary between $C$ and $D$. Consequently, for each
covering $\mathcal{C}$ of $%
\mathbb{R}
^{2}$ by a set of complete folding curves which have the same type and
define the same set $E_{1}$, there exists a unique fonction $\sigma _{%
\mathcal{C}}:%
\mathbb{Z}
^{2}\rightarrow \{-1,+1\}$ such that $\sigma _{\mathcal{C}}(X)=\sigma
_{C}(X) $\ for each curve $C\in \mathcal{C}$ and each vertex $X$ of $C$%
.\bigskip

\noindent \textbf{Lemma 3.4.} Let$\ \mathcal{C}$ and $\mathcal{D}$ be
coverings of $%
\mathbb{R}
^{2}$ by sets of complete folding curves. Suppose that all the curves have
the same type, define the same sets $E_{k}$ and have locally isomorphic
associated sequences. Then, for each $n\in 
\mathbb{N}
^{\ast }$, each $X\in 
\mathbb{Z}
^{2}-E_{2n-1}$ and each $U\in 2^{n}%
\mathbb{Z}
^{2}$, we have $\sigma _{\mathcal{C}}(X)=\sigma _{\mathcal{D}}(U+X)$.\bigskip

\noindent \textbf{Proof.}\ For each integer $k\leq 2n-1$, we have $%
F_{k}+2^{n}%
\mathbb{Z}
^{2}=F_{k}$ and therefore $U+F_{k}=F_{k}$. In particular, $X$ and $U+X$
belong to $F_{m}$ for the same integer $m\leq 2n-2$. As each curve of $%
\mathcal{C}$ and each curve of $\mathcal{D}$ have the same type and define
the same $E_{1}$, the map $Z\rightarrow U+Z$ induces a bijection from the
set of all supports of segments of $\mathcal{C}$ to the set of all supports
of segments of $\mathcal{D}$ which respects the orientation of the segments.

We consider a curve $C=(C_{i})_{i\in 
\mathbb{Z}
}\in \mathcal{C}$ indexed in such a way that the terminal point of $%
C_{2^{m}} $ is $X$, and a curve $D=(D_{i})_{i\in 
\mathbb{Z}
}\in \mathcal{D}$ indexed in such a way that the support of $D_{2^{m}}$\ is $%
U+S$, where $S$ is the support of $C_{2^{m}}$. The initial point $Y$ of $%
C_{1}$ and the initial point $Z$ of $D_{1}$ belong to $E_{m+1}$ since $X$
and $U+X$\ belong to $F_{m}$.

The sequences $(\zeta _{i})_{i\in 
\mathbb{Z}
}$ and $(\eta _{i})_{i\in 
\mathbb{Z}
}$ associated to $C$ and $D$ are locally isomorphic. Consequently, we have $%
(\zeta _{1},...,\zeta _{2^{m}-1})=(\eta _{1},...,\eta _{2^{m}-1})$ since $X$
and $U+X$\ belong to $F_{m}$, and therefore $%
(D_{1},...,D_{2^{m}})=U+(C_{1},...,C_{2^{m}})$ and $Z=U+Y$. As $%
U+F_{m+1}=F_{m+1}$, it follows that $Y$ and $Z$ both belong to $F_{m+1}$, or
both belong to $E_{m+2}$. In both cases, we have $(\zeta _{1},...,\zeta
_{2^{m+1}-1})=(\eta _{1},...,\eta _{2^{m+1}-1})$ since $C$ is locally
isomorphic to $D$, and therefore $\sigma _{\mathcal{C}}(X)=\zeta
_{2^{m}}=\eta _{2^{m}}=\sigma _{\mathcal{D}}(U+X)$.~~$\blacksquare $\bigskip

\noindent \textbf{Theorem 3.5.} Let $\mathcal{C}$ be a covering of $%
\mathbb{R}
^{2}$ by a set of complete folding curves. Then $\mathcal{C}$ satisfies the
local isomorphism property if and only if all the curves have the same type,
define the same sets $E_{k}$ and have locally isomorphic associated
sequences.\bigskip

\noindent \textbf{Proof.} The condition is necessary according to
Proposition 3.3 and the remark after Theorem 3.2. Now we show that it is
sufficient.

For each $X=(x,y)\in 
\mathbb{Z}
^{2}$ and each $k\in 
\mathbb{N}
^{\ast }$, we write\ $S_{X,k}=[x,x+k]\times \lbrack y,y+k]$, and we denote
by $\mathcal{E}_{X,k}$ the set of all segments of $\mathcal{C}$ whose
supports are contained in $S_{X,k}$. It suffices to prove that, for each $%
X\in 
\mathbb{Z}
^{2}$ and each $k\in 
\mathbb{N}
^{\ast }$, there exists $l\in 
\mathbb{N}
^{\ast }$ such that each $S_{Y,l}$\ contains some $Z\in 
\mathbb{Z}
^{2}$ with $\mathcal{E}_{Z,k}=(Z-X)+\mathcal{E}_{X,k}$.

We consider the largest integer $m$ such that $S_{X,k}$\ contains a point of 
$F_{m}$, and an integer $n$ such that $m\leq 2n-2$. For each $U\in 2^{n}%
\mathbb{Z}
^{2}$, we have, according to Lemma 3.4, $\sigma _{\mathcal{C}}(U+Y)=\sigma _{%
\mathcal{C}}(Y)$ for each $Y\in 
\mathbb{Z}
^{2}-E_{2n-1}$, and in particular for each $Y\in S_{X,k}$ which does not
belong to $E_{\infty }$.

If $S_{X,k}\cap E_{\infty }=\emptyset $, it follows that $\mathcal{E}%
_{U+X,k}=U+\mathcal{E}_{X,k}$\ for each $U\in 2^{n}%
\mathbb{Z}
^{2}$, since the curves of $\mathcal{C}$ have the same type. Then we have
the required property for $l=2^{n}$.

If $S_{X,k}$ contains the unique point $W$ of $E_{\infty }$, then we still
have $\mathcal{E}_{U+X,k}=U+\mathcal{E}_{X,k}$\ for each $U\in 2^{n}%
\mathbb{Z}
^{2}$\ such that $\sigma _{\mathcal{C}}(U+W)=\sigma _{\mathcal{C}}(W)$,\
since the curves of $\mathcal{C}$ have the same type. Moreover, we have $%
E_{2n}=W+2^{n}%
\mathbb{Z}
^{2}$\ since $W$ belongs to $E_{2n}$. We consider $V\in 2^{n}%
\mathbb{Z}
^{2}$ such that $V+W\in F_{2n}$ and $\sigma _{\mathcal{C}}(V+W)=\sigma _{%
\mathcal{C}}(W)$.\ We have $V+W+2^{n+1}%
\mathbb{Z}
^{2}=\{Y\in 
\mathbb{Z}
^{2}\mid Y\in F_{2n}$ and $\sigma _{\mathcal{C}}(Y)=\sigma _{\mathcal{C}%
}(W)\}$, and therefore $\mathcal{E}_{U+V+X,k}=U+V+\mathcal{E}_{X,k}$\ for
each $U\in 2^{n+1}%
\mathbb{Z}
^{2}$. Consequently we have the required property for $l=2^{n+1}$.~~$%
\blacksquare $\bigskip

Now we consider the coverings which consist of one complete folding curve.
The following result is a particular case of Theorem 3.5. It applies to the
folding curves considered in Theorem 3.7 and Example 3.8 below.\bigskip

\noindent \textbf{Corollary 3.6.} Any covering of $%
\mathbb{R}
^{2}$ by a complete folding curve satisfies the local isomorphism
property.\bigskip

\noindent \textbf{Theorem 3.7.} For each complete folding sequence $S$,
there exists a covering of $%
\mathbb{R}
^{2}$ by a complete folding curve associated to a sequence which is locally
isomorphic to $S$.\bigskip

\noindent \textbf{Proof.} Consider a curve $C=(C_{i})_{i\in 
\mathbb{Z}
}$\ associated to $S$. Let $\Omega $ consist of the finite curves parallel
to subcurves of $C$, such that $(0,0)$\ is one of their vertices. For each $%
F\in \Omega $, denote by $N(F)$ the largest integer $n$ such that $F$
contains a covering of $[-n,+n]^{2}$.

For any $F,G\in \Omega $, write $F<G$ if $F\subset G$ and $N(F)<N(G)$. As $S$
satisfies the local isomorphism property, the union of any strictly
increasing sequence in $\Omega $\ is a covering of $%
\mathbb{R}
^{2}$ by a complete folding curve associated to a sequence which is locally
isomorphic to $S$.

It remains to be proved that, for each $F\in \Omega $, there exists $G>F$ in 
$\Omega $. According to Proposition 2.6, there is an integer $m$ such that
each $(C_{i+1},...,C_{i+m})$ contains a curve parallel to $F$. By Theorem
3.1, there exists a finite subcurve $K$ of $C$ which contains a covering of
a square $X+[-n,+n]^{2}$ with $X\in 
\mathbb{Z}
^{2}$ and $n=m+N(F)+1$.

Let $i$ be an integer such that $X$ is the terminal point of $C_{i}$.
Consider a curve $H$ parallel to $F$\ and contained in $%
(C_{i+1},...,C_{i+m}) $. Denote by $\tau $ the translation such that $\tau
(F)=H$. Then $\tau ^{-1}(X)$\ belongs to $[-m,+m]^{2}$ since $\tau ((0,0))$
belongs to $X+[-m,+m]^{2}$. For $G=\tau ^{-1}(K)$, we have $F\subset G$ and $%
G$ covers $\tau ^{-1}(X)+[-n,+n]^{2}$. In particular, $G$ covers $%
[-N(F)-1,+N(F)+1]^{2}$.~~$\blacksquare $\bigskip

\noindent \textbf{Remark. }There is no covering of $%
\mathbb{R}
^{2}$ by a curve associated to $(\overline{S},+1,S)$ or to $(\overline{S}%
,-1,S)$, where $S$ is an $\infty $-folding sequence. In fact, such a curve $%
C $ would contain $4$ segments having the point of $E_{\infty }(C)$ as an
endpoint, and $E_{\infty }(\overline{S},+1,S)$\ (resp. $E_{\infty }(%
\overline{S},-1,S)$) would contain $2$ integers.\bigskip

\bigskip 
\begin{center}
\includegraphics[scale=0.16]{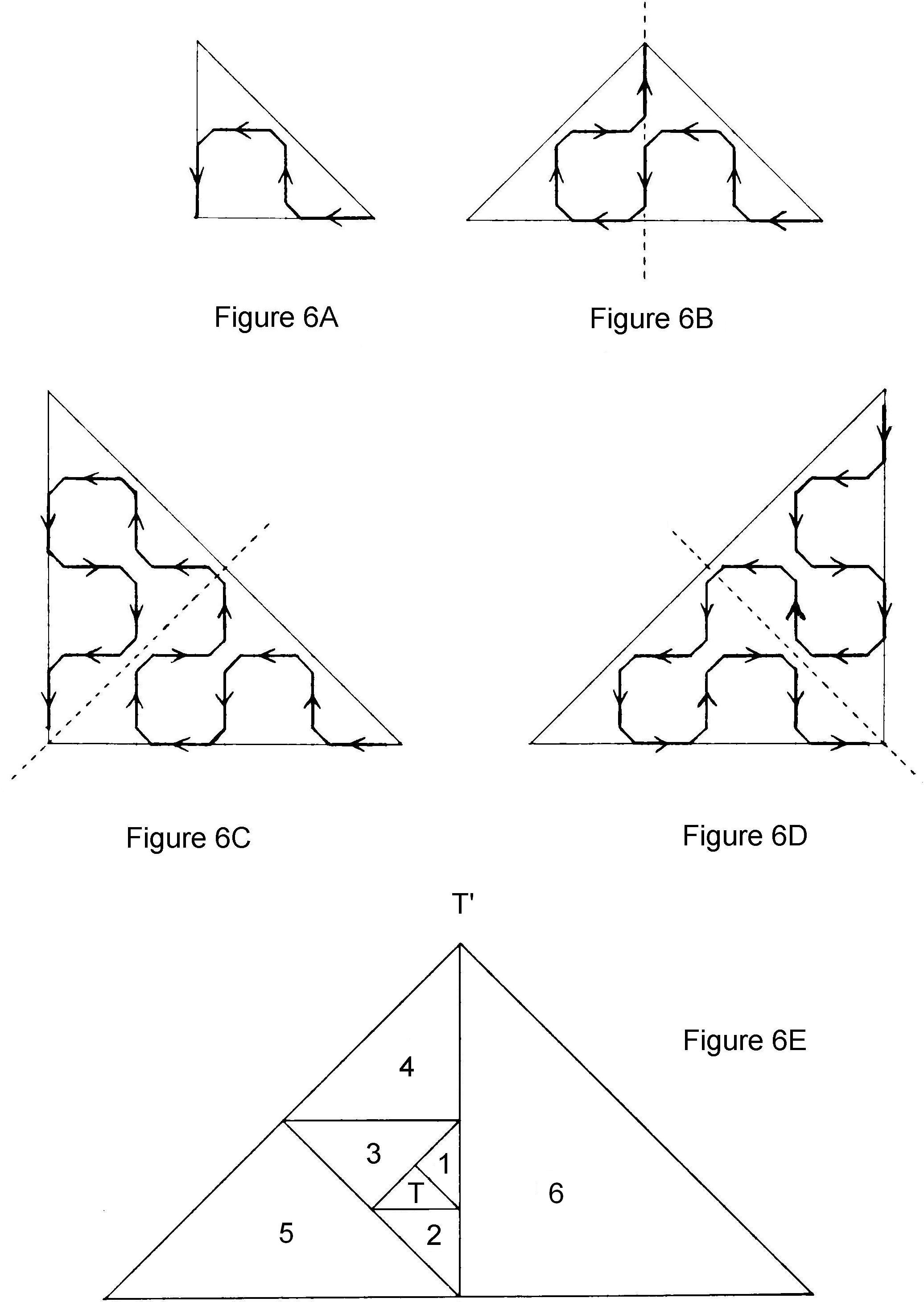}
\end{center}
\bigskip

\noindent \textbf{Example 3.8.}
There exists a covering of $%
\mathbb{R}
^{2}$ by a complete folding curve defined in an effective way.\bigskip

\noindent \textbf{Proof.} For each $(x,y)\in 
\mathbb{Z}
^{2}$ and each $n\in 
\mathbb{N}
^{\ast }$, we say that a $(2n)$-folding (resp. $(2n+1)$-folding) curve $C$
covers the isosceles right triangle $T=\langle
(x,y),(x+2^{n},y),(x,y+2^{n})\rangle $\ (resp. $T=\langle
(x,y),(x+2^{n},y+2^{n}),(x+2^{n+1},y)\rangle $) if it satisfies the
following conditions (cf. Fig. 6A, 6B, 6C):

\noindent - each interval $[(u,v),(u+1,v)]$ or $[(u,v),(u,v+1)]$ with $%
u,v\in 
\mathbb{Z}
$, contained in the interior of $T$, is the support of a segment of $C$;

\noindent - among the intervals $[(u,v),(u+1,v)]$ or $[(u,v),(u,v+1)]$ with $%
u,v\in 
\mathbb{Z}
$, contained in the same vertical or horizontal edge of $T$, alternatively
one over two is the support of a segment of $C$;

\noindent - no interval $[(u,v),(u+1,v)]$ or $[(u,v),(u,v+1)]$ with $u,v\in 
\mathbb{Z}
$, contained in the exterior of $T$, is the support of a segment of $C$;

\noindent - the vertex of the right angle of $T$ is the initial or the
terminal point of $C$;

\noindent - the vertex of one of the acute angles of $T$ is the initial or
the terminal point of $C$.

\noindent We extend this definition to the isosceles right triangles with
vertices in $%
\mathbb{Z}
^{2}$ obtained from $T$ by rotations of angles $\pi /2$, $\pi $, $3\pi /2$
(cf. Fig. 6D).

Now, let $T$ be one of the isosceles right triangles considered, let $k$ be
an integer, let $C$ be a $k$-folding curve which covers $T$, and let $S$ be
the sequence associated to $C$. If the initial (resp. terminal) point of $C$
is the vertex of the right angle of $T$, we associe to $(\overline{S},+1,S)$
and $(\overline{S},-1,S)$ (resp. $(S,+1,\overline{S})$ and $(S,-1,\overline{S%
})$) two $(k+1)$-folding curves $C_{1}$ and $C_{2}$ which contain $C$. In
both cases, the parts of $C_{1}$ and $C_{2}$ which correspond to $\overline{S%
}$ respectively cover the isosceles right triangles $T_{1}$ and $T_{2}$
which have one edge of their right angle in common with $T$, in such a way
that $C_{1}$ and $C_{2}$ respectively cover the isosceles right triangles $%
T\cup T_{1}$ and $T\cup T_{2}$ (cf. Fig. 6B, 6C and 6D).

For each $n\in 
\mathbb{N}
^{\ast }$ and each triangle $T=\langle
(x,y),(x+2^{n},y+2^{n}),(x+2^{n+1},y)\rangle $ with $(x,y)\in 
\mathbb{Z}
^{2}$, repete six times the operation above according to Figure 6E,
beginning with a $(2n+1)$-folding curve $C$ which covers $T$. Then we obtain
a $(2n+7)$-folding curve $C^{\prime }$ which contains of $C$. Moreover, $%
C^{\prime }$ covers a triangle $T^{\prime }=\langle (x^{\prime },y^{\prime
}),(x^{\prime }+2^{n+3},y^{\prime }+2^{n+3}),(x^{\prime }+2^{n+4},y^{\prime
})\rangle $ with $(x^{\prime },y^{\prime })\in 
\mathbb{Z}
^{2}$ and $T$ contained in the interior of $T^{\prime }$. By iterating this
process, we obtain a covering of $%
\mathbb{R}
^{2}$ by a complete folding curve.~~$\blacksquare $\bigskip

\noindent \textbf{Proposition 3.9.} Let$\ \mathcal{C}$ be a covering of $%
\mathbb{R}
^{2}$ by a set of complete folding curves which satisfies the local
isomorphism property. For each covering $\mathcal{D}$ of $%
\mathbb{R}
^{2}$, the following properties are equivalent:

\noindent 1) $\mathcal{D}$ consists of complete folding curves, and the
curves of $\mathcal{C}\cup \mathcal{D}$ have the same type, define the same
sets $E_{k}$ and have locally isomorphic associated sequences.

\noindent 2) $\mathcal{C}=\mathcal{D}$, or $E_{\infty }(\mathcal{C})\neq
\varnothing $ and $\mathcal{C},\mathcal{D}$ only differ in the way to
connect the four segments which have the unique point of $E_{\infty }(%
\mathcal{C})$ as an endpoint.\bigskip

\noindent \textbf{Proof.} If 1) is true, then 2) is also true since Lemma
3.4 implies $\sigma _{\mathcal{C}}(X)=\sigma _{\mathcal{D}}(X)$\ for each $%
X\in 
\mathbb{Z}
^{2}-E_{\infty }$. Conversely, if 2) is true, then 1) follows form the
remark after Corollary 1.9.~~$\blacksquare $\bigskip

\noindent \textbf{Remark.} If 1) and 2) are true, then\textbf{\ }$\mathcal{D}
$ satisfies the local isomorphism property by Theorem 3.5, and $\mathcal{D}$
is locally isomorphic to $\mathcal{C}$ by Theorem 3.2.\bigskip

\noindent \textbf{Theorem 3.10.} For each complete folding curve $C$, if $%
E_{\infty }(C)$ is empty or if the unique point of $E_{\infty }(C)$ is a
vertex of $C$, then $C$ is contained in a unique covering of $%
\mathbb{R}
^{2}$ by a set of complete folding curves which satisfies the local
isomorphism property. Otherwise, $C$ is contained in exactly two such
coverings, which only differ in the way to connect the four segments having
the unique point of $E_{\infty }$ as an endpoint.\bigskip

\noindent \textbf{Proof.} It suffices to show that $C$ is contained in a
covering of $%
\mathbb{R}
^{2}$ by a set of complete folding curves which satisfies the local
isomorphism property. In fact, for any such coverings $\mathcal{C},\mathcal{D%
}$, Proposition 3.3 and the remark after Theorem 3.2 imply that each curve
of $\mathcal{C}$ and each curve of $\mathcal{D}$ have the same type, define
the same sets $E_{k}$ and have locally isomorphic associated sequences. Then
the Theorem is a consequence of Proposition 3.9 and the Remark just after.

For each $m\in 
\mathbb{N}
^{\ast }$, denote by $\mathcal{C}_{m}$ the set of all segments of $C$ with
supports in $[-m,+m]^{2}$. Let $\Omega $ consist of the pairs $(\mathcal{E}%
,m)$, where $m\in 
\mathbb{N}
^{\ast }$ and $\mathcal{E}$ is a covering of $[-m,+m]^{2}$ containing $%
\mathcal{C}_{m}$, for which there exists $X\in 
\mathbb{Z}
^{2}$\ such that $X+\mathcal{E}\subset C$.

For any $(\mathcal{E},m),(\mathcal{F},n)\in \Omega $, write $(\mathcal{E}%
,m)<(\mathcal{F},n)$ if $\mathcal{E}\subset \mathcal{F}$ and $m<n$. If $(%
\mathcal{F},n)\in \Omega $ and $m\in \{1,...,n-1\}$, then we have $(\mathcal{%
E},m)\in \Omega $ and $(\mathcal{E},m)<(\mathcal{F},n)$ for the set $%
\mathcal{E}$ of all segments of $\mathcal{F}$ with supports in $[-m,+m]^{2}$.

First we show that, for each $m\in 
\mathbb{N}
^{\ast }$, there exists $(\mathcal{E},m)\in \Omega $. We can take $m$ large
enough so that $C$\ contains some segments with supports in $[-m,+m]^{2}$.
We consider a finite curve $A\subset C$\ which contains all these segments.
According to Proposition 2.6, there exists an integer $k$ such that each
subcurve of length $k$ of $C$ contains a curve parallel to $A$.

By Theorem 3.1, $C$ contains a covering of a square $X+[-k-2m,+k+2m]^{2}$
with $X\in 
\mathbb{Z}
^{2}$. The covering of $X+[-k,+k]^{2}$\ extracted from $C$ contains a curve
of length $\geq k$, which itself contains a curve $B$ parallel to $A$. We
consider $Y\in 
\mathbb{Z}
^{2}$\ such that $Y+A=B$. We have $Y\in X+[-k-m,+k+m]^{2}$, because $A$\
contains a point of $[-m,+m]^{2}$ and $B$ is contained in\ $X+[-k,+k]^{2}$.
Consequently, $C$ contains a covering $\mathcal{F}$ of $Y+[-m,+m]^{2}$. We
have $(\mathcal{E},m)\in \Omega $ for $\mathcal{E}=-Y+\mathcal{F}$ since $%
\mathcal{F}$ contains the set $Y+\mathcal{C}_{m}$ of all segments of $B$
with supports in $Y+[-m,+m]^{2}$.

Now, according to K\"{o}nig's Lemma, $\Omega $ contains a strictly
increasing sequence $(\mathcal{E}_{m},m)_{m\in 
\mathbb{N}
^{\ast }}$. The union $\mathcal{C}$ of the sets $\mathcal{E}_{m}$\ is a
covering of $%
\mathbb{R}
^{2}$ which contains $C$. Any finite curve $A\subset \mathcal{C}$ is
parallel to a subcurve of $C$ since it is contained in one of the sets $%
\mathcal{E}_{m}$. In particular, $\mathcal{C}$ is a covering of $%
\mathbb{R}
^{2}$ by a set of complete folding curves. It remains to be proved that $%
\mathcal{C}$ satisfies the local isomorphism property.

It suffices to show that, for each $m\in 
\mathbb{N}
^{\ast }$, there exists $n\in 
\mathbb{N}
^{\ast }$ such that each square $X+[-n,+n]^{2}$\ contains the image of $%
\mathcal{E}_{m}$\ under a translation. We consider a point $Y\in 
\mathbb{Z}
^{2}$\ such that $Y+\mathcal{E}_{m}\subset C$, and a finite curve $A\subset
C $\ which contains $Y+\mathcal{E}_{m}$. By Proposition 2.6, there exists $%
n\in 
\mathbb{N}
^{\ast }$\ such that each subcurve of length $n$ of $C$, and therefore each
subcurve of length $n$ of a curve of $\mathcal{E}$, contains a curve
parallel to $A$.\ Each square $X+[-n,+n]^{2}$\ contains a subcurve of length 
$n$ of a curve of $\mathcal{E}$, and therefore contains a curve parallel to $%
A$ and the image of $\mathcal{E}_{m}$\ under a translation.~~$\blacksquare $%
\bigskip

\noindent \textbf{Remark.} It follows that any $\infty $-folding curve $C$
is contained in exactly two coverings of $%
\mathbb{R}
^{2}$ by sets of complete folding curves which satisfy the local isomorphism
property. These coverings only differ in the way to connect the four
segments whose supports contain the origin of $C$.\bigskip

\noindent \textbf{Remark.} Let $C$ be a complete folding curve and let $S$
be the associated sequence. By Theorem 3.10, $C$ is contained in a covering $%
\mathcal{C}$ of $%
\mathbb{R}
^{2}$ by a set of complete folding curves which satisfies the local
isomorphism property. According to Theorem 3.7, there is also a covering $%
\mathcal{D}$ of $%
\mathbb{R}
^{2}$ by a complete folding curve $D$ associated to a sequence which is
locally isomorphic to $S$. Moreover, $\mathcal{D}$\ satisfies the local
isomorphism property by Corollary 3.6. If we choose $D$ with the same type
as $C$, then $\mathcal{C}$ and $\mathcal{D}$ are locally isomorphic
according to Theorem 3.2. On the other hand, $\mathcal{C}$ and $\mathcal{D}$
do not contain the same number of curves if $\{C\}$\ is not a covering of $%
\mathbb{R}
^{2}$.\bigskip

\noindent \textbf{Corollary 3.11.} 1) There exist $2^{\omega }$ pairwise not
locally isomorphic coverings of $%
\mathbb{R}
^{2}$ by sets of complete folding curves which satisfy the local isomorphism
property.

\noindent 2) If $\mathcal{C}$ is such a covering, then there exist $%
2^{\omega }$ isomorphism classes of coverings of $%
\mathbb{R}
^{2}$ which are locally isomorphic to $\mathcal{C}$.\bigskip

\noindent \textbf{Proof.} For any complete folding sequences $S,T$,\
consider two curves $C,D$\ associated to $S,T$ which have the same type. By
Theorem 3.10, there exist two coverings $\mathcal{C},\mathcal{D}$ of $%
\mathbb{R}
^{2}$, respectively containing $C,D$, by sets of complete folding curves
which satisfy the local isomorphism property.

By Theorem 3.2, $\mathcal{C}$ and $\mathcal{D}$ are locally isomorphic if
and only if $S$ and $T$ are locally isomorphic. In particular, the property
1) above follows from the property 1) of Theorem 1.11.

Moreover, if $\mathcal{C}$ and $\mathcal{D}$ are isomorphic, then $C$ is
isomorphic to one of the curves of $\mathcal{D}$,\ and $S$\ is isomorphic to
one of their associated sequences. Consequently, the property 2) above
follows from the property 2) of Theorem 1.11, since any covering of $%
\mathbb{R}
^{2}$ by a set of complete folding curves which satisfies the local
isomorphism property contains at most countably many curves (and in fact at
most $6$ curves by Theorem 3.12 below).~~$\blacksquare $\bigskip

Now we investigate the number of curves in a covering of $%
\mathbb{R}
^{2}$ by complete folding curves.\bigskip

\noindent \textbf{Theorem 3.12.} If a covering of $%
\mathbb{R}
^{2}$ by a set of complete folding curves satisfies the local isomorphism
property, then it contains at most $6$ curves.\bigskip

\noindent \textbf{Proof.} For each $X\in 
\mathbb{R}
^{2}$, each complete folding curve $C$ and each $n\in 
\mathbb{N}
$, write $\delta _{n}(X,C)=4^{-n}.\inf_{S\in E_{4n}(C)}d(X,S)$. We have $%
\delta _{n}(X,C)=\delta _{0}(X_{n},C^{(4n)})$, where $X_{n}$ is the image of 
$X$ in a representation of $C^{(4n)}$ which gives the length $1$ to the
supports of the segments.

Now consider $R\in E_{0}(C)$ and $S,T\in E_{4}(C)$ which are joined by a $4$%
-folding subcurve of $C$ having $R$ as a vertex. Then, according to Figure
3C, we have

\noindent $\inf (d(X,S),d(X,T))\leq \sqrt{(d(X,R)+3)^{2}+2^{2}}=\sqrt{%
d(X,R)^{2}+6.d(X,R)+13}$;

\noindent this maximum is reached with the second of the four $4$-folding
curves of Figure 3C, for each point $X$ which is at the left and at a
distance $\geq 3$ from the middle of the segment $ST$.

Consequently we have $\delta _{1}(X,C)\leq (1/4)\sqrt{\delta
_{0}(X,C)^{2}+6.\delta _{0}(X,C)+13}$. For $\delta _{0}(X,C)\leq 1,16$, it
follows

\noindent $\delta _{1}(X,C)\leq (1/4)\sqrt{1,16^{2}+6.1,16+13}<1,16$.

\noindent For $\delta _{0}(X,C)\geq 1,16$, it follows

\noindent $\delta _{1}(X,C)/\delta _{0}(X,C)\leq (1/4)\sqrt{1+6/\delta
_{0}(X,C)+13/\delta _{0}(X,C)^{2}}$

\noindent\ $\ \ \ \ \ \ \ \ \leq (1/4)\sqrt{1+6/1,16+13/1,16^{2}}<0,995$.

For each complete folding curve $C$, each $X\in 
\mathbb{R}
^{2}$ and each $n\in 
\mathbb{N}
$, the argument above applied to $C^{(n)}$ gives $\delta _{n+1}(X,C)<1,16$
if $\delta _{n}(X,C)<1,16$ and $\delta _{n+1}(X,C)<0,995.\delta _{n}(X,C)$
if $\delta _{n}(X,C)\geq 1,16$. In particular, we have $\delta
_{n}(X,C)<1,16 $ for $n$ large enough.

For each covering $\{C_{1},...,C_{k}\}$ of $%
\mathbb{R}
^{2}$ by a set of complete folding curves which satisfies the local
isomorphism property, consider $X\in 
\mathbb{R}
^{2}$ and $n\in 
\mathbb{N}
^{\ast }$ such that $\delta _{n}(X,C_{i})<1,16$ for $1\leq i\leq k$. Let $Y$
be the image of $X$ in a representation of $C_{1}^{(4n)},...,C_{k}^{(4n)}$
which gives the length $1$ to the supports of the segments. Then we have $%
\delta _{0}(Y,C_{i}^{(4n)})<1,16$ for $1\leq i\leq k$.

Write $Y=(y,z)$ and consider $u,v\in 
\mathbb{Z}
$ such that $\left\vert y-u\right\vert \leq 1/2$ and $\left\vert
z-v\right\vert \leq 1/2$. Then each $C_{i}^{(4n)}$ has a vertex among the
points listed below,\ since no other point $(w,x)\in 
\mathbb{Z}
^{2}$ satisfies $d((w,x),(y,z))<1,16$:

\noindent $(u,v),(u-1,v),(u+1,v),(u,v-1),(u,v+1)$ if $\left\vert
y-u\right\vert \leq 0,16$ and $\left\vert z-v\right\vert \leq 0,16$;

\noindent $(u-1,v-1),(u-1,v),(u-1,v+1),(u,v-1),(u,v),(u,v+1)$ if $y<u-0,16$;

\noindent $(u,v-1),(u,v),(u,v+1),(u+1,v-1),(u+1,v),(u+1,v+1)$ if $y>u+0,16$;

\noindent $(u-1,v-1),(u,v-1),(u+1,v-1),(u-1,v),(u,v),(u+1,v)$ if $z<v-0,16$;

\noindent $(u-1,v),(u,v),(u+1,v),(u-1,v+1),(u,v+1),(u+1,v+1),$ if $z>v+0,16$.

In the first case, there exist $12$ intervals of length $1$ with endpoints
in $%
\mathbb{Z}
^{2}$\ which have exactly one endpoint among the $5$ points considered. In
each of the four other cases, there exist $10$ intervals of length $1$ with
endpoints in $%
\mathbb{Z}
^{2}$\ which have exactly one endpoint among the $6$ points considered.

In all cases, each $C_{i}^{(4n)}$ has at least $2$ segments with supports
among these intervals. It follows $k\leq 6$ in the first case and $k\leq 5$
in each of the four other cases since, by Proposition 3.3, $%
C_{1}^{(4n)},...,C_{k}^{(4n)}$ are disjoint.~~$\blacksquare $\bigskip

\bigskip 
\begin{center}
\includegraphics[scale=0.16]{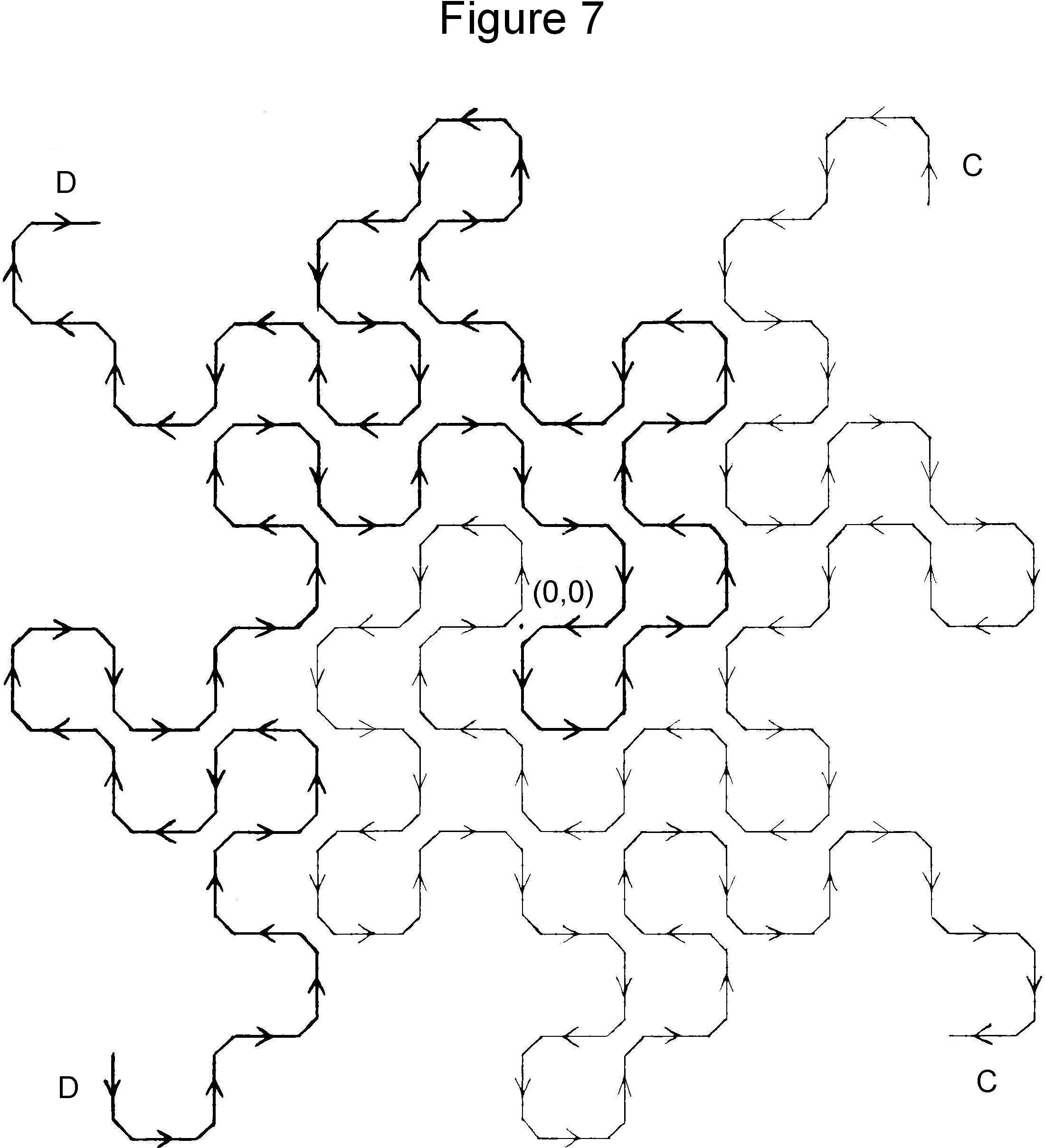}
\end{center}
\bigskip

\noindent \textbf{Example 3.13
(curves associated to the positive folding sequence).} The\emph{\ positive} 
\emph{folding sequence} mentioned in [4, p. 192] is the $\infty $-folding
sequence $R$ obtained as the direct limit of the $n$-folding sequences $%
R_{n} $ defined with $R_{n+1}=(R_{n},+1,\overline{R_{n}})$ for each $n\in 
\mathbb{N}
$. According to [3, Th. 4, p. 78], or by Theorem 3.15 below, there exists a
covering $\mathcal{C}$ of $%
\mathbb{R}
^{2}$ by $2$ complete folding curves $C,D$ associated to $S=(\overline{R}%
,+1,R)$ and having the same type (see Fig. 7). We have $E_{\infty
}(C)=E_{\infty }(D)=\{(0,0)\}$, and therefore $E_{k}(C)=E_{k}(D)$ for each $%
k\in 
\mathbb{N}
$. It follows from Theorem 3.5 that $\mathcal{C}$ satisfies the local
isomorphism property.\bigskip

\bigskip 
\begin{center}
\includegraphics[scale=0.14]{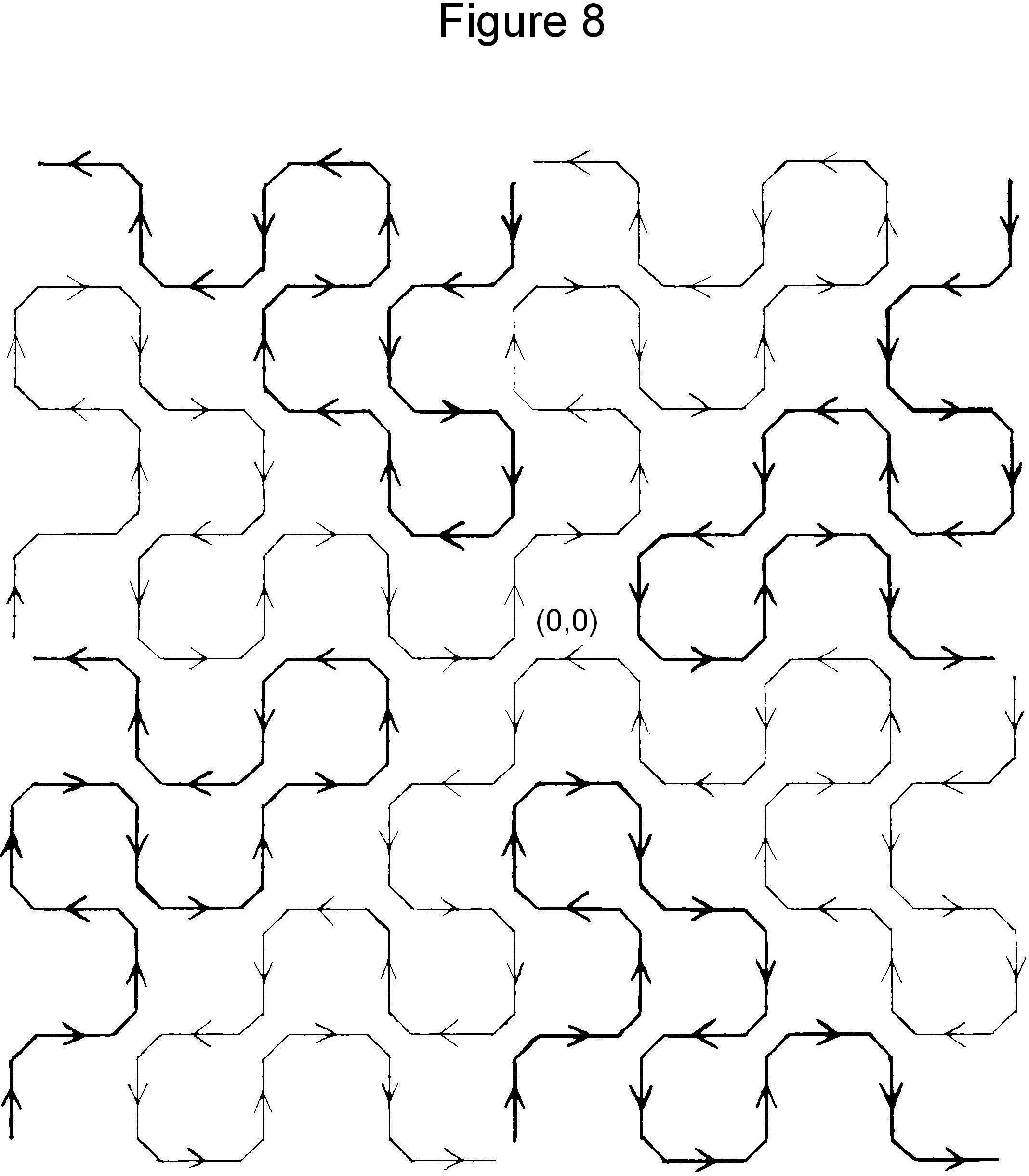}
\end{center}
\bigskip

\noindent \textbf{Exemple 3.14
(curves associated to the alternating folding sequence).} The\emph{\
alternating} \emph{folding sequence} described in [4, p. 134] is the $\infty 
$-folding sequence $R$ obtained as the direct limit of the $n$-folding
sequences $R_{n}$ defined with $R_{n+1}=(R_{n},(-1)^{n+1},\overline{R_{n}})$
for each $n\in 
\mathbb{N}
$. Contrary to Example 3.13, there exists no covering of $%
\mathbb{R}
^{2}$ by $2$ complete folding curves associated to $S=(\overline{R},+1,R)$.

Now, let $T$ be the complete folding sequence obtained as the direct limit
of the sequences $R_{2n}$, where each $R_{2n}$\ is identified to its second
copy in $R_{2n+2}=(R_{2n},-1,\overline{R_{2n}},+1,R_{2n},+1,\overline{R_{2n}}%
)$. Then there exists (cf. Fig. 8) a covering $\mathcal{C}$ of $%
\mathbb{R}
^{2}$ by $2$ curves associated to $S$, $2$ curves associated to $T$ and $2$
curves associated to $\overline{T}$, with all the curves having the same
type.

The $6$ curves are associated to locally isomorphic sequences. The point $%
(0,0)$\ belongs to $E_{\infty }$ in the $2$ curves which contain it. For
each $n\in 
\mathbb{N}
$, and in each of the $2$ curves which contain it, the point $(2^{n},0)$
(resp. $(-2^{n},0)$, $(0,2^{n})$, $(0,-2^{n})$) belongs to $F_{2n}$, while
the point $(2^{n},2^{n})$ (resp. $(-2^{n},2^{n})$, $(2^{n},-2^{n})$, $%
(-2^{n},-2^{n})$) belongs to $F_{2n+1}$. Consequently, the $6$ curves define
the same sets $E_{k}$, and $\mathcal{C}$ satisfies the local isomorphism
property by Theorem 3.5.\bigskip

We do not know presently if a covering $\mathcal{C}$ of $%
\mathbb{R}
^{2}$ by a set of complete folding curves which satisfies the local
isomorphism property can consist of $3$, $4$ or $5$ curves. If $E_{\infty }(%
\mathcal{C})\neq \varnothing $, then $\mathcal{C}$ consists of $2$ or $6$
curves according to the Theorem below:\bigskip

\noindent \textbf{Theorem 3.15.} Let $C$ be a curve associated to $S=(%
\overline{R},+1,R)$ or $S=(\overline{R},-1,R)$ where $R=(a_{h})_{h\in 
\mathbb{N}
^{\ast }}$\ is an $\infty $-folding sequence. Consider the unique covering $%
\mathcal{C}\supset C$ of $%
\mathbb{R}
^{2}$ by a set of complete folding curves which satisfies the local
isomorphism property.

\noindent 1) If $\{n\in 
\mathbb{N}
\mid a_{2^{n}}=(-1)^{n}\}$ is finite or cofinite, then $\mathcal{C}$
consists of $2$ curves associated to $S$ and $4$ other curves.

\noindent 2)\ Otherwise, $\mathcal{C}$ consists of $2$ curves associated to $%
S$.\bigskip

\noindent \textbf{Proof.} As the other cases can be treated in the same way,
we suppose that $C=(C_{h})_{h\in 
\mathbb{Z}
}$\ is a curve associated to $S=(\overline{R},+1,R)=(a_{h})_{h\in 
\mathbb{Z}
}$, and that $E_{\infty }(C)=\{(0,0)\}$. Then $\mathcal{C}$\ contains a
curve $D\neq C$\ such that $(0,0)$\ is a vertex of $D$. We write $%
D=(D_{h})_{h\in 
\mathbb{Z}
}$\ with $(0,0)$\ terminal point of $D_{0}$ and initial point of $D_{1}$.
The curve $D$ is associated to $S$, since it is associated to a sequence $%
T=(b_{h})_{h\in 
\mathbb{Z}
}$ with $T$ locally isomorphic to $S$, $E_{\infty }(T)=\{0\}$ and $b_{0}=+1$.

If $\{n\in 
\mathbb{N}
\mid a_{2^{n}}=(-1)^{n}\}$ is cofinite (resp. finite), we consider an odd
(resp. even) integer $k$ such that $a_{2^{n}}=(-1)^{n}$ (resp. $%
a_{2^{n}}=(-1)^{n+1}$) for $n\geq k$. In order to prove that $\mathcal{C}$\
contains $6$\ curves, it suffices to show that $C^{(k)}$ is contained in a
covering of $%
\mathbb{R}
^{2}$ which satisfies the local isomorphism property and which consists of $%
6 $\ complete folding curves. Consequently, it suffices to consider the case
where $a_{2^{n}}=(-1)^{n+1}$ for each $n\in 
\mathbb{N}
$. But this case is treated in Example 3.14 (see Fig. 8).

The proof of 2) uses arguments similar to those in the proof of Theorem 3.1.
For each $k\in 
\mathbb{N}
^{\ast }$, we say that a set $\mathcal{E}$ of disjoint curves covers $%
L_{k}=\{(u,v)\in 
\mathbb{R}
^{2}\mid \left\vert u\right\vert +\left\vert v\right\vert \leq k\}$ if each
interval $[(u,v),(u+1,v)]$ or $[(u,v),(u,v+1)]$ with $u,v\in 
\mathbb{Z}
$, contained in $L_{k}$, is the support of a segment of one of the curves.

For each set $\mathcal{E}$ of disjoint complete folding curves which have
the same type, define the same sets $E_{n}$ and are associated to locally
isomorphic sequences, if $\mathcal{E}^{(2)}$ covers $L_{k}$ for an integer $%
k\geq 2$, then $\mathcal{E}$ covers $L_{2k-1}$, and therefore covers $%
L_{k+1} $. By induction, it follows that, if $\mathcal{E}^{(2k)}$ covers $%
L_{2}$ for an integer $k\geq 2$, then $\mathcal{E}$ covers $L_{k+2}$.

For each $k\in 
\mathbb{N}
$, consider $l\geq 2k$ such that $a_{2^{l}}=a_{2^{l+1}}$. Then $%
\{(C_{-3}^{(l)},...,C_{4}^{(l)}),$\allowbreak $%
(D_{-3}^{(l)},...,D_{4}^{(l)})\}$ covers\ $L_{2}$\ (see Fig. 7 for the case $%
a_{2^{l}}=a_{2^{l+1}}=+1$). Consequently, $\{C^{(l-2k)},D^{(l-2k)}\}$ covers 
$L_{k+2}$, and the same property is true for $\{C,D\}$.~~$\blacksquare $%
\bigskip

\bigskip 
\begin{center}
\includegraphics[scale=0.13]{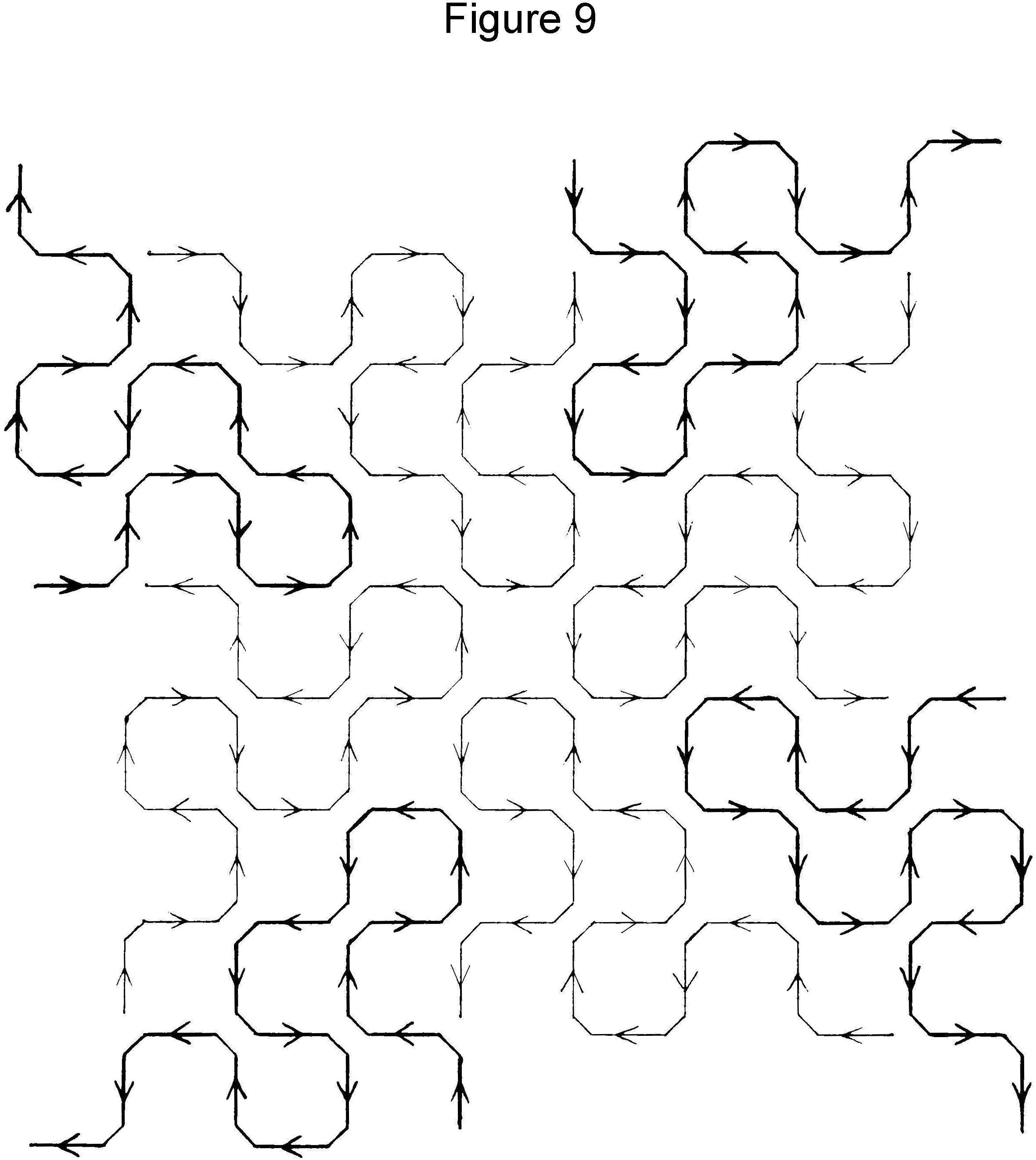}
\end{center}
\bigskip

A covering of $%
\mathbb{R}
^{2}$ by a set of complete folding curves can contain more than $6$ curves
if it does not satisfy the local isomorphism property. For the complete
folding sequence $T$ of Example 3.14, Figure 9 gives an example of a
covering of $%
\mathbb{R}
^{2}$ by $4$ curves associated to $T$ and $4$ curves associated to $-%
\overline{T}$ ($T$ and $-\overline{T}$ are not locally isomorphic by
Corollary 1.9). Anyway, the Theorem below implies that the number of
complete folding curves in a covering of $%
\mathbb{R}
^{2}$ is at most $24$:\bigskip

\noindent \textbf{Theorem 3.16.} There exist no more than $24$ disjoint
complete folding curves in $%
\mathbb{R}
^{2}$.\bigskip

In the proof of Theorem 3.16, we write $\delta ((x,y),(x^{\prime },y^{\prime
}))=\sup (\left\vert x^{\prime }-x\right\vert ,\left\vert y^{\prime
}-y\right\vert )$ for any $(x,y),(x^{\prime },y^{\prime })\in 
\mathbb{R}
^{2}$.\bigskip

\noindent \textbf{Lemma 3.16.1.} We have $\delta (X,Y)\leq 7.2^{n-2}-2$ for
each integer $n\geq 2$ and for any vertices $X,Y$ of a $(2n)$-folding
curve.\bigskip

\noindent \textbf{Proof of the Lemma.} The Lemma is true for $n=2$ according
to Figure 3C. Now we show that, if it is true for $n\geq 2$, then it is true
for $n+1$.

Let $C$ be a $(2n+2)$-folding curve, and let $Z_{0},...,Z_{2^{2n+2}}$ be its
vertices taken consecutively. Then $(Z_{4i})_{0\leq i\leq 2^{2n}}$ is the
sequence of all vertices of the $(2n)$-folding curve $C^{(2)}$ represented
on the same figure as $C$. It follows from the induction hypothesis applied
to $C^{(2)}$ that we have $\delta (X,Y)\leq 2(7.2^{n-2}-2)$ for any $X,Y\in
(Z_{4i})_{0\leq i\leq 2^{2n}}$. Moreover, for each vertex\ $U$ of $C$, there
exists $V\in (Z_{4i})_{0\leq i\leq 2^{2n}}$ such that $\delta (U,V)\leq 1$.
Consequently, we have $\delta (X,Y)\leq 2(7.2^{n-2}-2)+2=7.2^{n-1}-2$ for
any vertices $X,Y$ of $C$.~~$\blacksquare $\bigskip

\noindent \textbf{Lemma 3.16.2.} Let $n\geq 2$ be an integer and let $C$ be
a finite folding curve with two vertices $X,Y$ such that $\delta (X,Y)\geq
7.2^{n-2}-1$. Then $C$\ contains at least $2^{2n-1}$\ segments.\bigskip

\noindent \textbf{Proof of the Lemma.} By Lemma 3.16.1, we have $k\geq 2n+1$
for the smallest integer $k$ such that $C$ is contained in a $k$-folding
curve $D$. We consider the $(k-2)$-folding curves $D_{1},D_{2},D_{3},D_{4}$\
such that $D=(D_{1},D_{2},D_{3},D_{4})$. The curve $C$ contains one of the
curves\ $D_{2},D_{3}$ since the $(k-1)$-folding curves\ $(D_{1},D_{2})$, $%
(D_{2},D_{3})$, $(D_{3},D_{4})$ do not contain $C$. Consequently, $C$
contains a $(2n-1)$-folding curve, which consists of $2^{2n-1}$\ segments.~~$%
\blacksquare $\bigskip

\noindent \textbf{Proof of the Theorem.} Let $r\geq 2$ be an integer and let 
$C_{1},...,C_{r}$\ be disjoint complete folding curves. Consider an integer $%
k$ and some vertices $X_{1},...,X_{r}$ of $C_{1},...,C_{r}$\ belonging to $%
\left[ -k,+k\right] ^{2}$.

Now consider an integer $n\geq 2$ and write $N=7.2^{n-2}+k$. For each $i\in
\{1,...,r\}$, there exist some vertices $Y_{i},Z_{i}$ of $C_{i}$, with $%
\delta (0,Y_{i})=N$ and $\delta (0,Z_{i})=N$, such that $\left] -N,+N\right[
^{2}$ contains a subcurve of $C_{i}$ which has $X_{i}$ as a vertex and $%
Y_{i},Z_{i}$ as endpoints.

For each $i\in \{1,...,r\}$, we have $\delta (X_{i},Y_{i})\geq N-k=7.2^{n-2}$
and $\delta (X_{i},Z_{i})\geq N-k=7.2^{n-2}$. By Lemma 3.16.2, the part of
the subcurve of $C_{i}$\ between $X_{i}$ and $Y_{i}$ (resp. between $X_{i}$
and $Z_{i}$) contains at least $2^{2n-1}$\ segments.

As $C_{1},...,C_{r}$\ are disjoint, it follows that $\left] -N,+N\right[
^{2} $ contains at least $2r.2^{2n-1}=2^{2n}r$\ supports of segments of $%
C_{1}\cup ...\cup C_{r}$. But $\left] -N,+N\right[ ^{2}$ only contains $%
2(2N)(2N-1)<8N^{2}$ intervals of the form $](u,v),(u+1,v)[$ or $%
](u,v),(u,v+1)[$ with $u,v\in 
\mathbb{Z}
$. Consequently, we have $2^{2n}r<8N^{2}$ and $%
r<N^{2}/2^{2n-3}=(7.2^{n-2}+k)^{2}/2^{2n-3}$.

As the last inequality is true for each $n\geq 2$, it follows $r\leq 49/2<25$%
.~~$\blacksquare $\bigskip

\bigskip

\bigskip

\begin{center}
\textbf{References}
\end{center}

\bigskip

\noindent [1] J.-P. Allouche, The number of factors in a paperfolding
sequence, Bull. Austral. Math. Soc. 46 (1992), 23-32.

\noindent [2] A. Chang, T. Zhang, The fractal geometry of the
boundary of dragon curves, J. Recreational Math. 30 (1999), 9-22.

\noindent [3] C. Davis and D.E. Knuth, Number representations and
dragon curves I, J. Recreational Math. 3 (1970), 66-81.

\noindent [4] M. Dekking, M. Mend\`{e}s France, A. van der Poorten,
Folds!, Math. Intelligencer 4 (1982), 130-195.

\noindent [5] M. Mend\`{e}s France, A.J. van der Poorten, Arithmetic
and analytic properties of paperfolding sequences, Bull. Austral. Math. Soc.
24 (1981), 123-131.

\noindent [6] S.-M. Ngai, N. Nguyen, The Heighway dragon revisited,
Discrete Comput. Geom. 29 (2003), 603-623.

\noindent [7] F. Oger, Algebraic and model-theoretic properties of
tilings, Theoret. Comput. Sci. 319 (2004), 103-126.\bigskip
\bigskip

Francis Oger

C.N.R.S. - Equipe de Logique

D\'{e}partement de Math\'{e}matiques

Universit\'{e} Paris 7

case 7012, site Chevaleret

75205 Paris Cedex 13

France

E-mail: oger@logique.jussieu.fr.

\vfill

\end{document}